\documentclass{amsproc}
\usepackage{amssymb}
\usepackage[matrix, arrow, curve]{xy}
\usepackage{enumerate}
\newtheorem{theorem}{Theorem}[section]

\newtheorem{proposition}[theorem]{Proposition}

\newtheorem{definition}[theorem]{Definition}
\usepackage[pdftex]{graphicx,color}
\theoremstyle{definition}
\newtheorem{rem}[theorem]{Remark}

\newcommand{\F}{\mathcal F}
\newcommand{\E}{\mathcal E}
\newcommand{\D}{\mathcal D}
\newcommand{\K}{\mathcal K}

\theoremstyle{remark}

\numberwithin{equation}{section}

\begin{document}

\title[ $q$-Gamma weights and Painlev\'e]{Orthogonal polynomials on the unit circle,   $q$-Gamma weights, and
discrete Painlev\'e equations.}

\author{Philippe Biane}
\address{CNRS, IGM,
    Universit\'e Paris-Est, Champs-sur-Marne, FRANCE
}
\email{biane@univ-mlv.fr}
\thanks{This paper was written with the support of ANR GrandMa}
\subjclass{Primary ; Secondary }
\date{}
\begin{abstract}
We consider orthogonal polynomials on the unit circle with respect to a weight which is a quotient of $q$-gamma functions. We show that the Verblunsky coefficients of these polynomials satisfy  discrete Painlev\'e equations, in a Lax form, which correspond to an   $A_3^{(1)}$ surface in Sakai's classification.
\end{abstract}

\maketitle
\section{Introduction} 
 Motivated by number theoretic considerations, J.F. Burnol \cite{Bur02}, 
\cite{Bur04},
 has raised the problem of realizing the Fourier transform as a scattering. More precisely, consider the Fourier transform restricted to even functions on ${\bf R}$ which  is  given,  after appropriate normalization, by 
$${\mathcal F} f(x)=2\int_0^\infty \cos(2\pi xy)f(y)dy.$$
One can rewrite this as an isometry of the space 
$L^2(]0,+\infty[,dx/x)$, by putting $f(x)= g(x)/\sqrt x$, to obtain
$$\tilde{\mathcal F} g(x)=2\int_0^\infty \sqrt{xy}\cos(2\pi xy)g(y)dy/y$$
The interaction between the additive Fourier transform and the Mellin transform is captured by the formula, where the Fourier transform acts on tempered distributions,
\begin{equation}\label{four}
\tilde{\mathcal F} x^{+ik}=\pi^{-ik}\frac{\Gamma(1/4-ik/2)}{\Gamma(1/4+ik/2)}x^{-ik}\qquad k\ \text{real}
\end{equation}
which yields the gamma factor coming from the place at infinity in the functional equation for Riemann zeta function. Gamma factors for other $L$-functions are obtained by similar considerations involving the radial Fourier transform in 
${\bf R}^n$.

Interpreting $x^{ik}$ as an incoming wave, and $x^{-ik}$ as the corresponding outgoing  wave (in logarithmic coordinates) on the space $]0,+\infty[$,
it is a natural question to ask what is the scattering potential which realizes
 the map (\ref{four}) (we shall recall the basics of scattering theory in the next section).
This  was  answered  by Burnol \cite{BurP} who showed that, more generally, the function 
\begin{equation}\label{Gamma}
\frac{\Gamma(a-ik/2)\Gamma(b+ik/2)}{\Gamma(a+ik/2)\Gamma(b-ik/2)}
\end{equation} ($a,b$ real, $a>b$), is the scattering phase shift for a Dirac equation on $]0,+\infty[$, where the potential is given in terms of Painlev\'e transcendents.
In this paper we consider a discrete analogue of this problem. Discrete scattering problems have been introduced by Case and Kac \cite{Case-Kac} and studied further by Case and Geronimo \cite{Case-Geronimo}. In this last paper it is shown that, given a function $f$,  analytic inside the unit disk, continuous on the unit circle, then 
$f(e^{i\theta})/ \overline{f(e^{i\theta})}$
  is  the phase shift for a discrete scattering problem, and solving the inverse problem for this scattering problem 
 is equivalent to finding the recursion relations for the orthogonal polynomials on the unit circle with respect to the measure $d\theta/|f(e^{i\theta})|^2$. 
Recall that, if $\phi_n$ is the sequence of these polynomials, and 
$$\phi_n^*(z)=z^n\bar\phi_n(1/z)$$ then the recursion relation, due to Szeg\"o, can be written, with
$$\Phi_n(z)=\begin{pmatrix}\phi_n(z)
		\\
		\phi_n^*(z)\end{pmatrix}
$$
as
\begin{equation}\label{recintro}
\Phi_{n+1}(z)=B_n(z)\Phi_n(z)
\end{equation}
where
\begin{equation}\label{scat_disc}
B_n(z)=\begin{pmatrix}z&\alpha_{n+1}\\
z\bar\alpha_{n+1}&1\end{pmatrix}.
\end{equation}
The $\alpha_n$ are the Verblunsky coefficients of the sequence $\phi_n$,
see e.g. \cite{Sz},  \cite{Geronimus}, or \cite{Simon1}, \cite{Simon2}.
Equation (\ref{recintro}) is a discrete analogue of a Dirac system leading to scattering. 
Thus a natural problem, in view of 
 Burnol's  question, is to find the recursion coefficients for a weight $w$, 
expressed in terms of $q$-Gamma functions, of the form 
$$
w(z)=\frac{(az,q)_{\infty}( \bar az^{-1},q)_{\infty}}{(bz,q)_{\infty}
 ( \bar bz^{-1},q)_{\infty}}.
$$
with, as usual,
$$(z,q)_\infty=\prod_0^\infty(1-zq^i).$$
Indeed $w(e^{i\theta})=1/|f(e^{i\theta})|^2$ with $f(e^{i\theta})=(be^{i\theta},q)_{\infty}/(ae^{i\theta},q)_{\infty}$, and taking a suitable limit of $f(e^{i\theta})/ \overline{f(e^{i\theta})}$ for $q\to 1$ yields the phase shift (\ref{Gamma}).

 The weight $w$ satisfies a $q$-difference equation
$$w(qz)=\frac{V(z)}{W(z)}w(z)$$
where $V$ and $W$ are two polynomials of degree 2.
It turns out that a method, which goes back to Laguerre, further used by Freud \cite{Fr}, allows, in this situation, to obtain nonlinear relations between these recursion coefficients. This kind of question has seen renewed interest recently, see e.g. \cite{Magnus}, \cite{F-W},
\cite{vanassche}.
We shall  prove that
$$
\Phi_{n}(qz)=\frac{1}{V(z)}A_n(z)\Phi_n(z)
$$
where $A_n$ is a matrix with polynomial entries, of degree 2, which can be computed from the coefficients of the orthogonal polynomials. Compatibility between this equation and (\ref{recintro})  implies
the relation
$$
A_{n+1}(z)B_n(z)=B_n(qz)A_n(z).
$$ This last 
equation can be interpreted as  an isomonodromy  deformation of the $q$ difference operator associated with the matrix $A_n$. 
It is similar,  although it is not equivalent to,
 the Lax form of the discrete Painlev\'e
 VI equation studied by Jimbo and Sakai in \cite{JS}. We shall study it in details
and  prove that the equation that we obtain in this way is a Lax form for a discrete Painlev\'e 
equation, corresponding to the $A_3^{(1)}$ surface, according to the classification of Sakai
\cite{S}. As we shall see, this equation looks  like the one of \cite{JS}, but is obtained from  translation in a different direction in the affine Coxeter group $D_5^{(1)}$. 

This paper is organized as follows. We start by recalling some basic material on one-dimensional scattering theory, following \cite{F}, for the continuous case, and \cite{Case-Geronimo} for the discrete case. We also recall some of the theory of orthogonal polynomials on the unit circle, following closely \cite{F-W}.
 Laguerre's method is used in section 3 to derive nonlinear equations satisfied by the recursion coefficients of these orthogonal polynomials. Then, in section 4, we
give an explicit formula for these nonlinear relations as a birational map in two variables, and  show that  they give rise to  a regular isomorphism between two surfaces obtained by blowing up ${\bf P}^1\times {\bf P}^1$ at eight points. In section 5 we recall Sakai's theory, and  use it to identify our nonlinear relation as a 
  discrete Painlev\'e equation. We determine its type in Sakai's classification, and find its expression as a certain translation in the affine Weyl group of type 
$D_5^{(1)}$. Finally in the last section we derive a differential system by taking a continuous limit of the discrete equations.

The results of this paper are based on numerous lengthy computations. In order not to bore the reader we have omitted the details of most of these computations, but given enough details so that they can be reproduced.

I would like to thank Jean-Fran\c cois Burnol for keeping me informed about his work
which motivated me to undertake the study of these equations. I thank also
Alexei Borodin and Philip Boalch for their useful comments at an early stage of this work.
\section{One-dimensional scattering, continuous and discrete}
\subsection{Scattering for the Schr\"odinger equation}
One starts with the  equation
\begin{equation}\label{Schr}
-\frac{d^2}{dx^2}\psi(x,s)+q(x)\psi(x,s)=s^2\psi(x,s)\qquad s\in{\bf C}
\end{equation}
on the half line $[0,+\infty[$. This equation has a unique solution, $\psi(x,s)$, satisfying the boundary condition 
\begin{equation}
\psi(0,s)=0, \quad \frac{d}{dx}\psi(0,s)=1.
\end{equation}
For a potential $q$ decreasing sufficiently at infinity 
there are also two solutions, 
the Jost solutions $f(x,\pm s)$, determined by the conditions
\begin{eqnarray}
f(x,\pm s)e^{\mp isx}\to_{x\to\infty} 1.
\end{eqnarray}
For real $k$, 
the Jost solutions $f(x,k)$ and $f(x,-k)$  correspond to incoming and outgoing waves  and $\psi$ is  a linear combination
\begin{equation}\label{jost}
\psi(x,k)=\frac{1}{2ik}\left(f(x,k)\bar M(k)-f(x,-k)M(k)\right)
\end{equation}
such that
$M$ is the boundary value of the function $M(s)=f(0,s)$, analytic in the upper half-plane. 
It follows that, as $x\to \infty$,
$$\psi(x,k)= \frac{|M(k)|}{k}\sin(kx-\eta(k))+o(1)$$
where  $\eta(k)=\arg M(k)$ is the phase shift created by the potential $q$, with respect to the case $q=0$. The problem of inverse scattering is, given the phase shift $\eta$, to reconstruct the potential $q$. For phase shifts satisfying a growth condition at infinity, there exists a unique solution without bound states, which can be recovered by the methods developped in the 50's, see e.g. \cite{GL} and \cite{F}.

\subsection{Discrete scattering}
We now follow the exposition of Case and Geronimo \cite{Case-Geronimo}.
More information on orthogonal polynomials on the unit circle can be found  e.g.
in  Szeg\"o \cite{Sz}, Geronimus \cite{Geronimus}, or the comprehensive treatise by Simon \cite{Simon1},    \cite{Simon2}.

 Let $\mu$ be a finite positive  measure on the unit circle $U(1)$. Denote its moments
$$\int_0^{2\pi} e^{-ik\theta}d\mu(\theta)=c_k$$
and  $\phi_n(z)$ the sequence of orthogonal polynomials obtained by the Gram-Schmidt procedure from 
 $1,z,z^2,\ldots$, normalized by 
 $$\phi_n(z)=z^n+\ldots$$
This normalization  will be more convenient for our purposes than the usual one, where the polynomials have $L^2$-norm 1.
Denote
$$\sigma_n=\int |\phi_n|^2d\mu\quad \text{and}\quad
\alpha_n=\phi_n(0).$$
The $\alpha_n$ are the Verblunsky  coefficients of the system, cf \cite{Simon1}.
They satisfy
$$\sigma_{n+1}=(1-|\alpha_{n+1}|^2)\sigma_n$$ 
thus the sequence $\sigma_n$ is determined by $\sigma_0=\mu(U(1))$ and the sequence $\alpha_n$.
The orthogonal polynomials are subject to the recursion relation
\begin{equation}\label{recphi}
\left\{
\begin{array}{rcl}
\phi_{n+1}(z)&=&
z\phi_n(z)+\alpha_{n+1}z^n\bar\phi_n(1/z)\\ 
\\ 
\bar \phi_{n+1}(1/z)&=&z^{-n}\phi_n(z)+
\bar\alpha_{n+1}z^{-1}\bar\phi_n(1/z)
\end{array}\right .
\end{equation}
In matrix form, putting
$$\Psi(z,n)=\begin{pmatrix}\phi_n(z)\\ \bar\phi_n(1/z)\end{pmatrix}$$
and 
$$D(n)=\begin{pmatrix}z&\alpha_{n+1}z^n\\\bar\alpha_{n+1}z^{-n}&z^{-1}\end{pmatrix}$$
one can write
$$\Psi(z,n+1)=D(n)\Psi(z,n)$$

Assuming that the potential $\alpha$ satisfies  $\lim_{n\to\infty}\alpha_n=0$ fast enough,  
 one can construct the Jost solutions to equation (\ref{recphi})
$$
\Psi_{+}(z,n)=\begin{pmatrix}\phi_+z,n)\\ \hat\phi_+(z,n)\end{pmatrix}
\qquad
\Psi_{-}(z,n)=\begin{pmatrix}\hat\phi_-(z,n)\\ \phi_-(z,n)\end{pmatrix}
$$
 satisfying the boundary conditions
\begin{eqnarray*}
\lim_n |\phi_{+}(z,n)-z^{ n}|=0&\quad |z|>1\\
 \lim_n |\phi_{-}(z,n)-z^{- n}|=0& \quad|z|<1\\
\lim_n |\hat\phi_{+}(z,n)|=0 &\quad|z|>1  \\
\lim_n |\hat \phi_{-}(z,n)|=0&\quad |z|<1.
\end{eqnarray*} 
The solutions $\Psi_{\pm}$ are linearly independent, and
\begin{equation}\label{jostdisc}
\Psi(z,n)=f_+(z)\Psi_+(z,n)+f_-(z)\Psi_-(z,n)
\end{equation}
where $$f_+(z)=\bar f_-(z) \quad\text{for}\quad |z|=1$$
Equation (\ref{jostdisc}) is the discrete analogue of (\ref{jost}).

The function  $f_+$ has analytic continuation inside the unit disc,
it  is related to the original measure $\mu$ on the unit circle by
$$d\mu(\theta)=\frac{d\theta}{|f_+(e^{i\theta})|^2}$$
The phase $\eta(\theta)=\arg(f_+(e^{i\theta}))$ is the phase shift induced by the potential $\alpha_n$. It can be recovered from the function
$\frac{1}{|f_+(e^{i\theta})|^2}$ on the unit circle.

One sees that reconstructing the phase shift from the potential $\alpha$ is equivalent to finding the orthogonality relation for the sequence of polynomials $\phi_n$, solution to the recursion equation (\ref{recphi}). Conversely, the inverse scattering problem is the easy problem of finding the recursion coefficients for the orthogonal polynomials, knowing their orthogonality measure. In fact the connection between scattering theory and orthogonal   polynomials has been known for a long time and is at the heart of the Gelfand-Levitan method (see \cite{GL}).
\subsection{Comparing the discrete and continuous cases}
Given a solution $A$ of 
the Riccatti equation
$$q=A'+A^2$$
the equation (\ref{Schr}), for $s=k$ real,
is equivalent to the Dirac system
\begin{equation}\label{Dir2}\left\{
\begin{array}{rcl}
dy/dx&=&Ay+kz,\\ 
-dz/dx&=&Az+ky.
\end{array}\right .
\end{equation}
The Dirac system (\ref{Dir2}) can then be put in the form
\begin{equation}\label{Wsys}\left\{
\begin{array}{rcl}
dW/dx=A\bar We^{2ikx},\\ 
d\bar W/dx=AWe^{-2ikx}.
\end{array}\right .
\end{equation}
with 
$$W=(y+iz)e^{ikx}$$
Similarly, putting
$$\omega_n(\theta)=e^{-in\theta}\phi_n(e^{i\theta})$$
 the system (\ref{recphi}) becomes
\begin{equation}
\label{Omegasys}\left\{
\begin{array}{rcl} 
\omega_{n+1}-\omega_n&=&
\alpha_{n+1}e^{-in\theta}\bar\omega_n,\\  
\bar \omega_{n+1}-\bar\omega_n&=&
\bar \alpha_{n+1}e^{in\theta}\omega_n.
\end{array}
\right .
\end{equation}
It is now clear that (\ref{Omegasys}) is a discretized version of 
(\ref{Wsys}).

\subsection{Recursion relations and the Caratheodory function}
We recall a few properties of the orthogonal polynomials which will be needed in the sequel.

The function 
$$F(z)=\int\frac{e^{i\theta}+z}{e^{i\theta}-z}
d\mu(\theta)$$
is called the Caratheodory function of $\mu$, it is holomorphic in ${\bf C}^*\setminus U(1)$, and 
has the expansion
$$\begin{array}{rcl}
F(z)&=&c_0+2\sum_{k=1}^\infty c_kz^k\quad |z|<1\\
F(z)&=&-c_0-2\sum_{k=1}^\infty c_{-k}z^{-k}\quad |z|>1
\end{array}$$
Let 
$$
\phi_n^*(z)=z^n\bar\phi_n(1/z)
$$
and
introduce the associated polynomials
\begin{eqnarray*}\psi_n(z)&=&\int\frac{e^{i\theta}+z}{e^{i\theta}-z}(\phi_n(e^{i\theta})-
\phi_n(z))
d\mu(\theta)\\
\psi_0(z)&=&1\\
\psi_n^*(z)&=&-\int\frac{e^{i\theta}+z}{e^{i\theta}-z}(z^n\overline{
\phi_n(e^{i\theta})}-
\phi_n^*(z))
d\mu(\theta)\\
\psi_0^*(z)&=&1
\end{eqnarray*}
and the functions
\begin{equation}\nonumber
\epsilon_n(z)=\psi_n(z)+F(z)\phi_n(z)=
\int\frac{e^{i\theta}+z}{e^{i\theta}-z}\phi_n(e^{i\theta})d\mu(\theta)
\end{equation}
\begin{equation}\nonumber
\epsilon^*_n(z)=\psi^*_n(z)-F(z)\phi^*_n(z)
=\int\frac{e^{i\theta}+z}{e^{i\theta}-z}\overline{\phi(e^{i\theta})}z^nd\mu(\theta)
\end{equation}
then, for $n\geq 1$,
\begin{eqnarray}\label{evaleps1}
\epsilon_n(z)&=&2\sigma_nz^n
+O(z^{n+1})\quad z\to 0
\\
\label{evaleps2}
\epsilon_n(z)&=&2\sigma_n\alpha_{n+1}z^{-1}+
O(z^{-2})\quad z\to \infty 
\\
\label{evaleps3}
\epsilon_n^*(z)&=&
2\sigma_n\bar\alpha_{n+1}z^{n+1}+
O(z^{n+2})\quad z\to 0
\\ \label{evaleps4}
\epsilon^*_n(z)&=&2\sigma_n+O(z^{-1})
\quad z\to \infty
\end{eqnarray}
In other words, the fractions 
$-\frac{\psi_n}{\phi_n}$ and
$ \frac{\psi^*_n}{\phi_n^*}$ 
are  reduced continued  fractions for $F$, near  $\infty$ and $0$
respectively. Furthermore,
$$
\begin{array}{rcl}
\epsilon_{n+1}(z)&=&
z\epsilon_n(z)-\alpha_{n+1}\epsilon_n^*(z)\\
\\
\epsilon_{n+1}^*(z)&=&
\epsilon_n^*(z)-
\bar\alpha_{n+1}z\epsilon_n(z).
\end{array}$$
The 
recursion  relations can be put in the form
\begin{equation}\label{recursion}Y_{n}(z)=\begin{pmatrix}\phi_{n}(z)&\frac{\epsilon_n(z)}{w(z)}\\
\phi_n^*(z)&-\frac{\epsilon_n^*(z)}{w(z)}\end{pmatrix}
\end{equation}
\begin{equation}\label{eqY}
Y_{n+1}(z)=B_n(z)Y_n(z).\end{equation}
with
\begin{equation}\label{B*}
B_n(z)=\begin{pmatrix}z&\alpha_{n+1}\\
z\bar\alpha_{n+1}&1\end{pmatrix}
\end{equation}

There are  also discrete Wronskian identities, derived from the recursion relations,  namely
\begin{eqnarray}\label{Wphipsi}
\phi_{n+1}(z)\psi_n(z)-\psi_{n+1}(z)\phi_n(z)&=&
\phi_{n+1}(z)\epsilon_n(z)-\epsilon_{n+1}(z)\phi_n(z)\\&=&{2\alpha_{n+1}}\nonumber
{\sigma_n}z^n
\end{eqnarray}
\begin{eqnarray}\label{Wphipsi*}
\phi^*_{n+1}(z)\psi^*_n(z)-\psi^*_{n+1}(z)\phi^*_n(z)&=&
\phi^*_{n+1}(z)\epsilon^*_n(z)-\epsilon^*_{n+1}(z)\phi^*_n(z)\\&=&\nonumber
2\bar\alpha_{n+1}
\sigma_nz^{n+1}
\end{eqnarray}
\begin{eqnarray}\label{Wphi*psi}
\phi_{n}(z)\psi^*_n(z)+\psi_{n}(z)\phi^*_n(z)
&=&\phi_{n}(z)\epsilon^*_n(z)+\epsilon_{n}(z)\phi^*_n(z)\\&=&\nonumber
2\sigma_nz^n.
\end{eqnarray}

\section{Nonlinear relations for the Verblunsky coefficients}
\subsection{Laguerre's method}
We now adapt the method of Laguerre \cite{L} and Freud \cite{Fr} for finding nonlinear relations between the Verblunsky coefficients. The computations that follow are inspired by  Forrester and Witte \cite{F-W}, who treat the case of difference equations 
of the form $w(z+1)=\rho(z)w(z)$.
Another derivation for the $q$-recurrence relations can be found in the paper    \cite{IW} by Ismail and Witte, however we shall need them in the form below, which is slightly different from the results of \cite{IW}, so we give a detailed derivation here.
 Let us  assume that the measure $\mu$ has the form
$$d\mu(\theta)=w(e^{i\theta})d\theta$$
 with $w$ holomorphic in ${\bf
C^*}$, satisfying
\begin{equation}\label{eqw}
w(qz)=\rho(z)w(z)
\end{equation}
where $q\in ]0,1[$,  and $\rho$ is a rational function,
\begin{equation}\rho(z)=\frac{V(z)}{W(z)}\end{equation}
 for some polynomials $V,W$.
The case of interest is when
 
\begin{equation}\label{w}
w(e^{i\theta})=\left|\frac{(ae^{i\theta},q)_{\infty}}{(be^{i\theta},q)_{\infty}}\right|^2=\frac{(ae^{i\theta},q)_{\infty}( \bar ae^{-i\theta},q)_{\infty}}{( be^{i\theta},q)_{\infty}
 (\bar be^{-i\theta},q)_{\infty}}
\end{equation}
 $q\in]0,1[,\, a,b\in {\bf C^*}$,
 thus
$$w(z)=\frac{(az,q)_{\infty}( \bar az^{-1},q)_{\infty}}{(bz,q)_{\infty}
 ( \bar bz^{-1},q)_{\infty}}$$
satisfies
 $$w(qz)=\frac{(qz-\bar a)(1- bz)}{(qz-\bar b)( az-1)}w(z),$$
 and we may take
 \begin{eqnarray}\label{VV}V(z)&=&(qz-\bar a)( bz-1)\\ \label{WW}
 W(z)&=&(qz-\bar b)(1- az).\end{eqnarray}

For the moment we shall make the computations with general polynomials $V$ and $W$, and come back to the specific values (\ref{VV}), (\ref{WW})
in the last section.

Equation (\ref{eqw}) implies that the
  Caratheodory function satisfies 

\begin{equation}\label{eqF}
W(z)F(qz)=V(z)F(z)+U(z)
\end{equation}
 with $U$ a polynomial of  degree at most $D=\max(\text{deg}\, W,\text{deg}\, V)$.
Indeed, take the partial fraction expansion  
$$\rho(z)=\delta(z)+\sum_j\frac{\omega_j}{z-a_j},\qquad \delta\ \text{polynomial of degree less than $(\text{deg}\, V-\text{deg}\, W)_+$}$$
with $W(z)=C\prod_j(z-a_j)$, then
$$
\begin{array}{rcl}
F(qz)&=&\int \frac{\zeta+qz}{\zeta-qz}w(\zeta)\frac{d\zeta}{i\zeta}\\
\\
&=&\int \frac{\zeta+z}{\zeta-z}w(q\zeta)\frac{d\zeta}{i\zeta}\qquad
\text{by a change of integration contour}
\\
\\
&=&\int \frac{\zeta+z}{\zeta-z}\rho(\zeta)w(\zeta)\frac{d\zeta}{i\zeta}\\
\\
&=&\rho(z)F(z)+\int \frac{\zeta+z}{\zeta-z}(\rho(\zeta)-\rho(z))
w(\zeta)\frac{d\zeta}{i\zeta}\\
\\
&=&\rho(z)F(z)+\int\left[\tilde\delta(\zeta,z)+\sum_j\frac{\omega_j}{z-a_j} \frac{\zeta+z}{\zeta-a_j}
w(\zeta)\right]\frac{d\zeta}{i\zeta}
\\ &&\qquad\qquad\qquad \tilde\delta(\zeta,z)=(\zeta+z)\frac{\delta(\zeta)-\delta(z)}{\zeta-z}
\\
\\
&=&(V(z)F(z)+U(z))/W(z)
\end{array}
$$
where $U$ is a   polynomial of degree at most
$ D$.
We shall now use (\ref{eqF}) and the characterization of the quotients
$\psi_n/\phi_n$ and $\psi_n^*/\phi_n^*$ as reduced continued fraction expansions.
Let us start from 
\begin{eqnarray*}
\epsilon_n(z)&=&\psi_n(z)+F(z)\phi_n(z)\\ \text{and}
\\  
\epsilon_n(qz)&=&\psi_n(qz)+(\frac{V(z)}{W(z)}F(z)+\frac{U(z)}{W(z)})\phi_n(qz)
\quad\text{by}\ (\ref{eqF})\end{eqnarray*}

and eliminate $F(z)$ between these two equations, to get
$$\begin{array}{c}
W(z)\epsilon_n(qz)\phi_n(z)-V(z)\epsilon_n(z)\phi_n(qz)=\\
W(z)\psi_n(qz)\phi_n(z)-V(z)\psi_n(z)\phi_n(qz)+U(z)\phi_n(z)\phi_n(qz).
\end{array}$$
The right hand side is a 
 polynomial in $z$ and, according to
(\ref{evaleps1})-(\ref{evaleps4}),
the left hand side is $O(z^n)$ near 0 and $O(z^{D+n-1})$ near $\infty$,
therefore there exists a polynomial $\Theta_n(z)$ of degree less than $D-1$, such that
\begin{equation}\label{theta}
2\alpha_{n+1}
\sigma_nz^n\Theta_n(z)=
W(z)\epsilon_n(qz)\phi_n(z)-V(z)\epsilon_n(z)\phi_n(qz).
\end{equation}
In  case $\alpha_{n+1}=0$, one has
$\phi_{n+1}(z)=z\phi_n(z)$ and $\phi^*_{n+1}(z)=\phi^*_n(z)$, and the sought relation for $\phi_{n+1}(qz)$ can be obtained from that of $\phi_{n}(qz)$. 
From 
$$2\alpha_{n+1}
\sigma_nz^n\Theta_n(z)=
W(z)\psi_n(qz)\phi_n(z)-V(z)\psi_n(z)\phi_n(qz)+U(z)\phi_n(z)\phi_n(qz)$$
we get, according to (\ref{Wphipsi}),
$$\begin{array}{rcl}[\phi_{n+1}(z)\psi_n(z)-\psi_{n+1}(z)\phi_n(z)]\Theta_n(z)
&=&W(z)\psi_n(qz)\phi_n(z)-V(z)\psi_n(z)\phi_n(qz)\\&&
+U(z)\phi_n(z)\phi_n(qz)\end{array}$$
or
$$\psi_n(z)(\phi_{n+1}(z)\Theta_n(z)+V(z)\phi_n(qz))
=\phi_n(z)(W(z)\psi_n(qz)+
U(z)\phi_n(qz)+\psi_{n+1}\Theta_n(z)).$$
This expression is a polynomial, common multiple of $\phi_n$ and $\psi_n$
which have no  common zero (see (\ref{Wphipsi})), so it can be put in the form
$$\Omega_n(z)\phi_n(z)\psi_n(z)$$ where $\Omega_n$ is a polynomial of degree at most
$ D$.
Thus $\Omega_n$ satisfies
\begin{eqnarray}\label{eqomega1}
\Omega_n(z)\phi_n(z)&=&\phi_{n+1}(z)\Theta_n(z)+V(z)\phi_n(qz)\\
\label{eqomega2}
\Omega_n(z)\psi_n(z)&=&W(z)\psi_n(qz)+
U(z)\phi_n(qz)+\psi_{n+1}(z)\Theta_n(z)
\end{eqnarray}
One obtains $\Omega_n$ by multiplying (\ref{eqomega2}) by $\phi_{n+1}$, 
substracting  (\ref{eqomega1}) multiplied by $\psi_{n+1}$, and using
(\ref{Wphipsi}), thus
\begin{equation}\label{defomega}
2\alpha_{n+1}
\sigma_nz^n\Omega_n(z)=W(z)\epsilon_n(qz)
\phi_{n+1}(z)-V(z)\epsilon_{n+1}(z)\phi_n(qz)
\end{equation}

Equation (\ref{eqomega1}) gives
\begin{equation}\label{recphi1V}
V(z)\phi_n(qz)=\Omega_n(z)\phi_n(z)-\phi_{n+1}(z)\Theta_n(z)
\end{equation}
or, using (\ref{recphi}),
\begin{equation}\label{recphi2V}
V(z)\phi_n(qz)=(\Omega_n(z)-z\Theta_n(z))\phi_n(z)
-
\alpha_{n+1}\Theta_n(z)\phi_n^*(z)
\end{equation}
On multiplying (\ref{eqomega1}) by $F(z)$ and adding (\ref{eqomega2})
 $$
\Omega_n(z)\epsilon_n(z)=
V(z)\epsilon_n(qz)\frac{w(z)}{w(qz)}+\epsilon_{n+1}(z)\Theta_n(z)
$$
therefore
$$
\Omega_n(z)\frac{\epsilon_n(z)}{w(z)}=
V(z)\frac{\epsilon_n(qz))}{w(qz)}+\frac{\epsilon_{n+1}(z)}{w(z)}\Theta_n(z).
$$
Finally 
\begin{equation}\label{recepsV}
V(z)\frac{\epsilon_n(qz))}{w(qz)}=\Omega_n(z)
\frac{\epsilon_n(z)}{w(z)}-\frac{\epsilon_{n+1}(z)}{w(z)}\Theta_n(z)
\end{equation}
and
\begin{equation}\label{epsilonq}
V(z)\frac{\epsilon_n(qz))}{w(qz)}=(\Omega_n(z)
-z
\Theta_n(z))\frac{\epsilon_n(z)}{w(z)}+
\alpha_{n+1}\Theta_n(z)\frac{\epsilon^*_{n}(z)}{w(z)}
\end{equation}

Equations for $\phi_n^*$ and $\epsilon_n^*$ are obtained by similar arguments.
One finds polynomials $\Theta^*_n(z)$ and $\Omega^*_n(z)$ with the same bounds on their degrees as  $\Theta_n(z)$ and $\Omega_n(z)$, such that

\begin{equation}\label{eqphi*}
V(z)\phi^*_n(qz)=(\Omega^*_n(z)-
\Theta^*_n(z))\phi^*_n(z)
-
z\bar\alpha_{n+1}\Theta^*_n(z)\phi_n(z)
\end{equation}
and
\begin{equation}\label{epsilonq*}
V(z)\frac{\epsilon^*_n(qz)}{w(qz)}=(\Omega^*_n(z)
-
\Theta^*_n(z))\frac{\epsilon^*_n(z)}{w(z)}+
z\bar\alpha_{n+1}\Theta^*_n(z)\frac{\epsilon_{n}(z)}{w(z)}
\end{equation}

\begin{proposition}
Let 
 \begin{equation}
\label{An}A_n(z)=\begin{pmatrix}
\Omega_n(z)-z\Theta_n(z)&
\alpha_{n+1}\Theta_n(z)
\\
-z\bar\alpha_{n+1}\Theta^*_n(z)
&\Omega^*_n(z)-\Theta^*_n(z)
\end{pmatrix}.
\end{equation}

then, with $Y_n$ defined by (\ref{recursion}),
the   relations (\ref{recphi2V}),
 (\ref{epsilonq}), (\ref{eqphi*}), (\ref{epsilonq*}) can be put in matrix form as

\begin{equation}
\label{eqYq}
Y_{n}(qz)=\frac{1}{V(z)}A_n(z)Y_n(z)
\end{equation}
Compatibility relations between 
 (\ref{eqY}) and (\ref{eqYq}) imply that

\begin{equation}\label{fund}
A_{n+1}(z)B_n(z)=B_n(qz)A_n(z).
\end{equation}
\end{proposition}

The coefficients of the polynomials $\Theta_n,\Theta_n^*,\Omega_n,\Omega_n^*$ in (\ref{An}) can be computed from the coefficients of the polynomials $\phi_n$.
For example, with $V$ and $W$ given by  (\ref{VV}) and (\ref{WW}) one finds, putting
 $$\Theta_n(z)=\lambda z+\mu$$

 \begin{equation}\label{eqtheta}
2\alpha_{n+1}z^n\Theta_n(z)=
W(z)\epsilon_n(qz)\phi_n(z)-V(z)\epsilon_n(z)\phi_n(qz)
\end{equation}
 and comparing  the two expressions near zero,
$$2\alpha_{n+1}z^n\mu=2\alpha_nz^n(
\bar bq^n-\bar a
)$$ 
Comparing at infinity gives $\lambda$, and
finally
\begin{equation}\label{theta_n}\Theta_n(z)=
(a
-bq^{n+1})z+(-\bar a+\bar bq^n)\frac{\alpha_n}{\alpha_{n+1}}.
\end{equation}
Similar reasoning yield the other polynomials in terms of the Verblunsky coefficients.
Using
$$\phi_n(z)=z^n+\beta_nz^{n-1}+\ldots+\alpha_n$$
with
$$\beta_n=\sum_{j=1}^n\alpha_j\bar\alpha_{j-1}$$
One finds
\begin{equation}\label{omega_n}
\Omega_n(z)-z\Theta_n(z)=bq^{n+1}z^2+tz+\bar bq^n
\end{equation}
where $t$ is a complicated expression in terms of $\alpha_n,\alpha_{n+1}$ and $\beta_n$. Similarly
\begin{eqnarray}\label{theta_n*}\Theta^*_n(z)=(aq-bq^{n+1})\frac{\bar\alpha_n}{\bar\alpha_{n+1}}z+bq^{n+1}-\bar a\\\label{omega_n*}
\Omega^*_n(z)-\Theta^*_n(z)=aqz^2+t^*z+\bar a
\end{eqnarray}
Equation (\ref{fund}) implies nonlinear recursion relations among these coefficients.
\subsection{The isomonodromy deformation}
Rather than do the tedious computations here, we shall, in the next section, 
start from the relation (\ref{fund}) and derive everything from its qualitative features. This will allow us to make the connection with Sakai's theory of discrete Painlev\'e equations. Before this we first make a few remarks.
Equation  (\ref{fund}) has the form
\begin{equation}\label{fund2}\tilde A(z)=B(qz)A(z)B(z)^{-1}\end{equation}
where $$A(z)=A_dz^d+\ldots +A_0,\quad\tilde A(z)=\tilde A_dz^d+\ldots +\tilde A_0$$
and $$B(z)=B\begin{pmatrix}z&0\\0&1\end{pmatrix}$$ for some constant matrix $B$.
Equation (\ref{fund2}) is an isomonodromy transformation of a  $q$-difference equation, as studied by
Birkhoff 
\cite{Bir}, and recently by Borodin \cite{Bor}. 
It is closely related to the equation appearing in \cite{JS}.
Such isomonodromy transformations, for difference equations instead of $q$-difference equations, have been studied in depth by Arinkin and Borodin in \cite{AB}, using geometric methods. These methods, involving the definition of ``$q$-connections'' should also apply to our situation, although we do not use them here.

In our case, the matrices $A_d,\tilde A_d$ are lower triangular matrices, and $A_0,\tilde A_0$ are upper triangular. Let $\kappa_1,\kappa_2$ be the eigenvalues of 
 $A_d$, and $\theta_1,\theta_2$ those of 
 $A_0, $ so that
$$A_d=\begin{pmatrix}\kappa_1&0\\ *&\kappa_2\end{pmatrix},\qquad 
A_0=\begin{pmatrix}\theta_1& *\\ 0&\theta_2\end{pmatrix}$$
Consider the matrix 
$$\hat A(z)=\begin{pmatrix}
qz& 0\\ 0&1 \end{pmatrix}A(z)\begin{pmatrix}
z^{-1}& 0\\ 0&1 \end{pmatrix}=\hat A_2z^2+\hat A_1z+\hat A_0$$
then 
$$\hat A_d=\begin{pmatrix}
q\kappa_1& *\\  0&
\kappa_2
\end{pmatrix},\qquad
\hat A_0=\begin{pmatrix}
q\theta_1& 0\\ *&
\theta_2
\end{pmatrix}, \qquad \text{det}\hat A(z)=q\text{det} A(z)$$

We can now conjugate $\hat A(z)$ by a constant matrix $B$ in order to get the matrix 
$\tilde A(z)$ with
$$\tilde A_d=\begin{pmatrix}q\kappa_1&0\\ *&\kappa_2\end{pmatrix},\qquad
\tilde A_0=\begin{pmatrix}q\theta_1& *\\ 0&\theta_2\end{pmatrix}.$$

The matrix $B$ is determined up to left multiplication by a diagonal invertible matrix. One can choose its coefficients to be expressed as rational functions of the coefficients of $A(z)$. Furthermore
 equation (\ref{fund2}) can be used in the same way  to deduce  $A(z)$ from $\tilde A(z)$.
It follows that   equation (\ref{fund2}) defines a birational transformation
between two spaces of matrices with polynomial coefficients. We shall give an explicit formula for this transformation in the next section.
The recursion equations for Verblunsky coefficients are given by iterating this transformation.
In the next section we shall investigate in details the case where $d=2$, which corresponds to the case of the Verblunsky coefficients for the  $q$-Gamma weight
 (\ref{w}).

\section{A birational transformation }

\subsection{The space $X_{\kappa_1,\kappa_2,\theta_1,\theta_2,c_1,c_2,c_3,c_4}$ }
We now embark on the analysis of equation (\ref{fund2}), in the case of degree 2 polynomial matrices. Consider a two-dimensional vector space $H$, and 
let $$A(z)=A_0+A_1z+A_2z^2$$ be an endomorphism of $H$ which depends polynomially on the variable $z$, with $A_0,A_1,A_2$, elements of $\text{End}(H)$. 
\begin{definition}
 We let 
$$\bar X_{\kappa_1,\kappa_2,\theta_1,\theta_2,c_1,c_2,c_3,c_4}$$ denote the space of all such $A(z)$ satisfying 
$$\{\kappa_1,\kappa_2\}=\text{sp}\,A_2,\qquad
\{\theta_1,\theta_2\}=\text{sp}\,A_0$$ and 
\begin{equation}\label{detA}
\text{det}\,A(z)=\kappa_1\kappa_2(z-c_1)(z-c_2)(z-c_3)(z-c_4)
\end{equation}
for some fixed complex constants  $\kappa_1,\kappa_2,\theta_1,\theta_2,c_1,\ldots
,c_4$, satisfying the equation 
\begin{equation}\kappa_1\kappa_2c_1c_2c_3c_4=\theta_1\theta_2\end{equation}
 (which comes from considering $\hbox{det}\,A(0)$).
We let  $$ X_{\kappa_1,\kappa_2,\theta_1,\theta_2,c_1,c_2,c_3,c_4}$$ 
denote its quotient  by the action of $GL(H)$  by conjugation.
These spaces are complex manifolds.
\end{definition}

In the sequel $c_1,\ldots,c_4$ will be fixed so we will omit them from the notation most of the time.

We consider the  case of generic parameters, in particular 
$\kappa_1\ne\kappa_2,\theta_1\ne\theta_2$, and  $A_0,A_2$ are semisimple. We denote 
$E_1,E_2$ the eigenspaces of $A_2$, and $F_1,F_2$ those of $A_0$.
The set of  $A(z)$ such that $E_2\ne F_1$ is a dense open set
$\Omega$, invariant by conjugation, and for each element of this open set we can chose a basis $(u,v)$  with  $u\in F_1,v\in E_2$ so that, in this basis,   
$A_0$ and $A_2$ have matrices of the form
$$A_2=\begin{pmatrix}
\kappa_1& 0\\  *&
\kappa_2
\end{pmatrix},\qquad
A_0=\begin{pmatrix}
\theta_1& *\\0&
\theta_2
\end{pmatrix}$$
As in the preceding section, let $q$ be a nonzero complex number and 
consider  the matrix 
$$\hat A(z)=\begin{pmatrix}
qz& 0\\ 0&1 \end{pmatrix}A(z)\begin{pmatrix}
z^{-1}& 0\\ 0&1 \end{pmatrix}=\hat A_2z^2+\hat A_1z+\hat A_0$$
then 
$$\hat A_2=\begin{pmatrix}
q\kappa_1& *\\  0&
\kappa_2
\end{pmatrix},\qquad
\hat A_0=\begin{pmatrix}
q\theta_1& 0\\ *&
\theta_2
\end{pmatrix}, \qquad \text{det}\hat A(z)=q\text{det} A(z)$$
therefore the matrix $\hat A(z)$ belongs to the space 
 $\bar X_{q\kappa_1,\kappa_2,q\theta_1,\theta_2}$.
Thus $A(z)\to\hat A(z)$ defines a map from a dense open subset of 
$\bar X_{\kappa_1,\kappa_2,\theta_1,\theta_2}$ to
$\bar X_{q\kappa_1,\kappa_2,q\theta_1,\theta_2}$.
Furthermore, as explained also at the end of the preceding section,  this map
 can be inverted, 
thus  defining a bijection between   dense open subsets of 
$\bar X_{\kappa_1,\kappa_2,\theta_1,\theta_2}$ and
$\bar X_{q\kappa_1,\kappa_2,q\theta_1,\theta_2}$. Clearly the map passes to the quotients and  induces a bijection
 between dense open subsets of $ X_{\kappa_1,\kappa_2,\theta_1,\theta_2}$ and
$ X_{q\kappa_1,\kappa_2,q\theta_1,\theta_2}$, which  is a birational map.
\begin{definition}
We call $\phi:X_{\kappa_1,\kappa_2,\theta_1,\theta_2}\to
 X_{q\kappa_1,\kappa_2,q\theta_1,\theta_2}$
the birational map constructed above.
\end{definition}
In the next section, we shall introduce a parametrization of the spaces in order to give an explicit formula for this birational transformation.
\subsection{Jimbo and Sakai  parametrization }
 We will use the following parametrization of the space 
$X_{\kappa_1,\kappa_2,\theta_1,\theta_2}$, after Jimbo and Sakai \cite{JS}
\begin{equation}\label{js}A(z)=\begin{pmatrix}
\kappa_1((z-y)(z-\alpha)+z_1)& z-y\\z(\gamma z+\beta)&
\kappa_2((z-y)(z-\delta)+z_2)
\end{pmatrix}
\end{equation}
Thus $y$ is the unique root of the polynomial which is the $12$ coefficient of $A$, and $\kappa_1z_1,\kappa_2z_2$ are the values of the $11$, resp. $22$ coefficient of $A(z)$ at $z=y$.
One has  $$y\alpha+z_1=\theta_1/\kappa_1,\qquad y\alpha+z_2=\theta_2/\kappa_2$$
The quantities  $y,z_1,z_2$ satisfy
\begin{equation}\label{kz}z_1z_2=(y-c_1)(y-c_2)(y-c_3)(y-c_4)\end{equation}
where the equality follows from considering the determinant of $A(y)$.
Let us  put $$\xi=(y-c_1)(y-c_2)/\kappa_1z_1=\kappa_2z_2/(y-c_3)(y-c_4).$$

 The quantities
 $\beta$ and $\gamma$ are determined as rational functions of $y,\xi$ using
(\ref{detA}).  

It follows that
the variables $(y,\xi)$ in ${\bf C}^*\times {\bf C^*}$  parametrize 
 a dense open subset of $X_{\kappa_1,\kappa_2,\theta_1,\theta_2}$. Using this parametrization, 
 the space $X_{\kappa_1,\kappa_2,\theta_1,\theta_2}$ is identified with  
${\bf P}^1\times {\bf P}^1$ blown up at the seven points  
$$\begin{matrix}(c_1,0),\quad (c_2,0),\quad (c_3,\infty),\quad (c_4,\infty)\\
(0,c_1c_2/\theta_1),\quad(0,c_1c_2/\theta_2),\quad(\infty,1/\kappa_1)
\end{matrix}$$ minus the strict transforms of the four lines 
$y=0,y=\infty$, $\xi=0,\xi=\infty$, see Fig. 2 below.

$$\resizebox{7 cm}{!}{\begin{picture}(0,0)%
\includegraphics{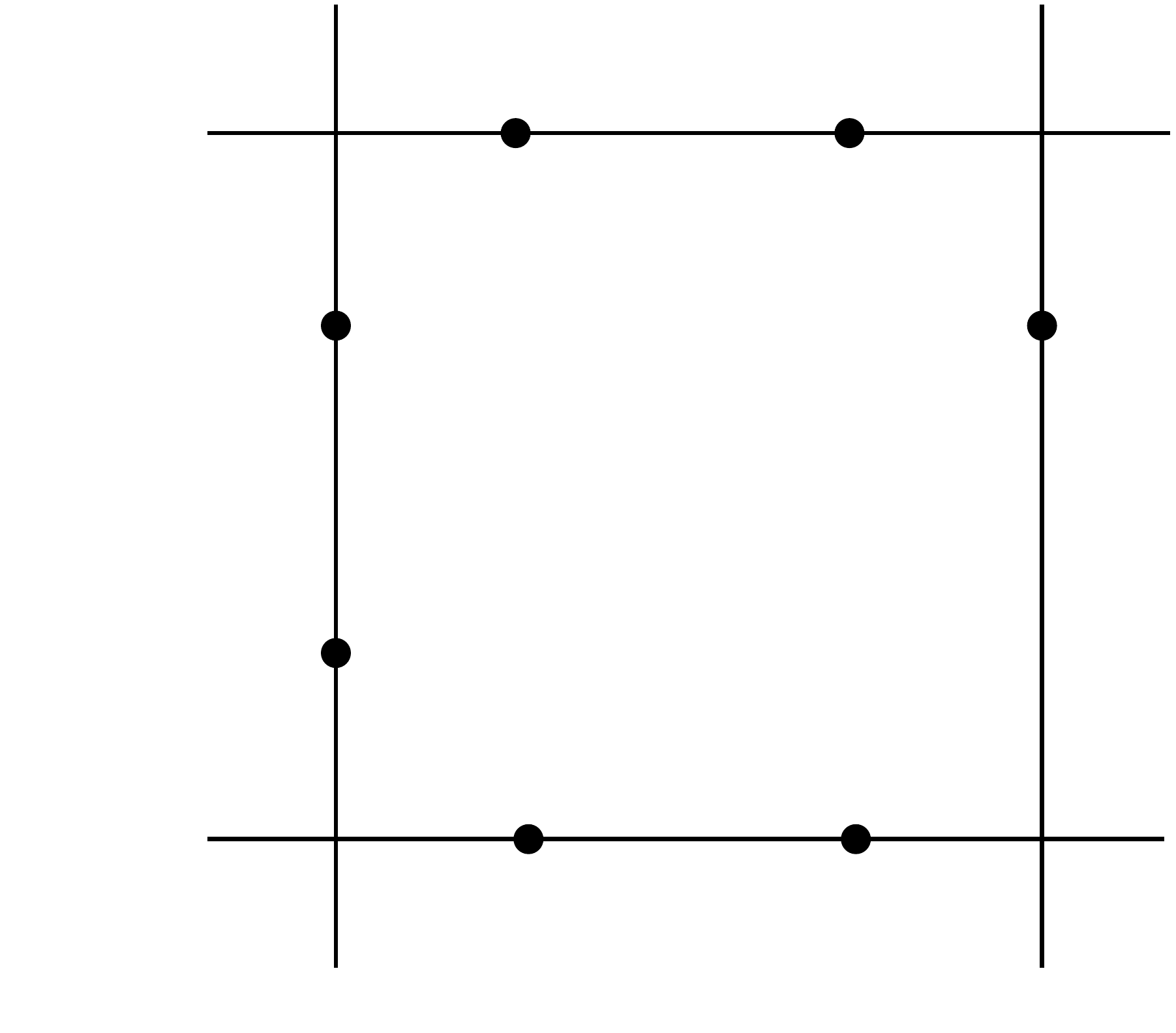}%
\end{picture}%
\setlength{\unitlength}{4144sp}%
\begingroup\makeatletter\ifx\SetFigFontNFSS\undefined%
\gdef\SetFigFontNFSS#1#2#3#4#5{%
  \reset@font\fontsize{#1}{#2pt}%
  \fontfamily{#3}\fontseries{#4}\fontshape{#5}%
  \selectfont}%
\fi\endgroup%
\begin{picture}(8238,7239)(796,-7267)
\put(4006,-6496){\makebox(0,0)[lb]{\smash{{\SetFigFontNFSS{20}{24.0}{\rmdefault}{\mddefault}{\updefault}{\color[rgb]{0,0,0}$(c_1,0)$}%
}}}}
\put(6256,-6541){\makebox(0,0)[lb]{\smash{{\SetFigFontNFSS{20}{24.0}{\rmdefault}{\mddefault}{\updefault}{\color[rgb]{0,0,0}$(c_2,0)$}%
}}}}
\put(3646,-646){\makebox(0,0)[lb]{\smash{{\SetFigFontNFSS{20}{24.0}{\rmdefault}{\mddefault}{\updefault}{\color[rgb]{0,0,0}$(c_3,\infty)$}%
}}}}
\put(6031,-646){\makebox(0,0)[lb]{\smash{{\SetFigFontNFSS{20}{24.0}{\rmdefault}{\mddefault}{\updefault}{\color[rgb]{0,0,0}$(c_4,\infty)$}%
}}}}
\put(1081,-2446){\makebox(0,0)[lb]{\smash{{\SetFigFontNFSS{20}{24.0}{\rmdefault}{\mddefault}{\updefault}{\color[rgb]{0,0,0}$(0,c_1c_2/\theta_1)$}%
}}}}
\put(1126,-4696){\makebox(0,0)[lb]{\smash{{\SetFigFontNFSS{20}{24.0}{\rmdefault}{\mddefault}{\updefault}{\color[rgb]{0,0,0}$(0,c_1c_2/\theta_2)$}%
}}}}
\put(8326,-2401){\makebox(0,0)[lb]{\smash{{\SetFigFontNFSS{20}{24.0}{\rmdefault}{\mddefault}{\updefault}{\color[rgb]{0,0,0}$(\infty,1/\kappa_1)$}%
}}}}
\put(2521,-7126){\makebox(0,0)[lb]{\smash{{\SetFigFontNFSS{20}{24.0}{\rmdefault}{\mddefault}{\updefault}{\color[rgb]{0,0,0}$y=0$}%
}}}}
\put(1216,-6046){\makebox(0,0)[lb]{\smash{{\SetFigFontNFSS{20}{24.0}{\rmdefault}{\mddefault}{\updefault}{\color[rgb]{0,0,0}$\xi=0$}%
}}}}
\put(811,-1051){\makebox(0,0)[lb]{\smash{{\SetFigFontNFSS{20}{24.0}{\rmdefault}{\mddefault}{\updefault}{\color[rgb]{0,0,0}$\xi=\infty$}%
}}}}
\put(7606,-7126){\makebox(0,0)[lb]{\smash{{\SetFigFontNFSS{20}{24.0}{\rmdefault}{\mddefault}{\updefault}{\color[rgb]{0,0,0}$y=\infty$}%
}}}}
\end{picture}%
} $$

\centerline{Fig. 2. The space $X_{\kappa_1}$}

\bigskip

The curves corresponding to
$F_2=E_1$ and $F_2=E_2$ are respectively  the strict transforms  of the points
$(0,c_1c_2/\theta_2)$ and $(0,c_1c_2/\theta_1)$.
\subsection{Explicit birational expression for $\phi$}We now describe the map $\phi$ as a birational transformation, using the Jimbo-Sakai coordinates (\ref{js}).
Let 
 $A(z)$ be as in (\ref{js}) and
consider the matrix
$$\begin{array}{rcl}\hat A(z)&=&\begin{pmatrix}
qz& 0\\ 0&1 \end{pmatrix}A(z)\begin{pmatrix}
z^{-1}& 0\\ 0&1 \end{pmatrix}\\&=&\begin{pmatrix}
q\kappa_1((z-y)(z-\alpha)+z_1)& qz(z-y)\\(\gamma z+\beta)&
\kappa_2((z-y)(z-\delta)+z_2)
\end{pmatrix}\end{array}$$
Then $\hat A(z)$ belongs to the space 
$X_{q\kappa_1,\kappa_2,q\theta_1,\theta_2}$. In order to find its Jimbo-Sakai coordinates, we have to use the basis formed by the spaces 
$\hat E_2$ and $\hat F_1$.
It follows that
 $$\tilde A(z)=B^{-1}\hat A(z)B$$ with
$$B=\begin{pmatrix}
q\theta_1-\theta_2
& -q\\ \beta&
q\kappa_1-\kappa_2
\end{pmatrix}$$
on the condition that the matrix $B$ so defined is invertible.
In order to go from  $\tilde A$ to $A$ one proceeds in the opposite way.
Start from 
$$\tilde A=\begin{pmatrix}
q\kappa_1((z-\tilde y)(z-\tilde\alpha)+\tilde z_1)& z-\tilde y\\\tilde z(\tilde\gamma z+\tilde\beta)&
\kappa_2((z-\tilde y)(z-\tilde\delta)+\tilde z_2)
\end{pmatrix}$$
conjugate it by a constant matrix, so that its    $21$
 coefficient has degree 1 and its  
$12$ coefficient is $z\times$(a polynomial with degree 1). For this use eigenvectors of  $\tilde A_2$ and $\tilde A_0$, thus
$\hat A(z)=C^{-1}\tilde A(z)C$
with
$$C=\begin{pmatrix}q\kappa_1-\kappa_2
& -{\tilde y}\\q {\tilde\gamma}&
q\theta_1-\theta_2
\end{pmatrix}$$
We will not compute the inverse transformation.

Let us now compute an explicit formula for $\phi$, using the coordinates $(y,\xi)$. We start from the matrix
$$A(z)=\begin{pmatrix}
\kappa_1z^2+az+\theta_1& z-y\\z(\gamma z+\beta)&
\kappa_2z^2+bz+\theta_2
\end{pmatrix}$$
with 

\begin{eqnarray}\label{a}a&=&-\kappa_1y+\frac{(y-c_1)(y-c_2)}{y\xi}-\frac{\theta_1}{y}\\
\label{b}
b&=&-\kappa_2y+\frac{\kappa_1\kappa_2\xi(y-c_3)(y-c_4)}{y}-\frac{\theta_2}{y}
\end{eqnarray}

and

\begin{eqnarray}
\label{beta}
\theta_1b+\theta_2a+\beta y&=&-\kappa_1\kappa_2\sigma_3
\\ \label{gamma}
\kappa_1b+\kappa_2a-\gamma &=&-\kappa_1\kappa_2\sigma_1
\end{eqnarray}
where
the $\sigma_i$ are the elementary symmetric functions in the $c_i$.
Then
\begin{equation}\label{atilde}\tilde A(z)=B(qz)A(z)B(z)^{-1}=\begin{pmatrix}
q\kappa_1z^2+\tilde az+q\theta_1& w(z-\tilde y)\\w^{-1}z(\tilde \gamma z+\tilde\beta)&
\kappa_2z^2+\tilde bz+\theta_2
\end{pmatrix}
\end{equation}
with
\begin{eqnarray*}
 B(z)&=&\begin{pmatrix}
z(q\kappa_1-\kappa_2)
&
q
\\
 -z\beta
&q\theta_1-\theta_2
\end{pmatrix}
\\
 B(z)^{-1}&=&\Delta^{-1}\begin{pmatrix}
z^{-1}(q\theta_1-\theta_2)
&
-z^{-1}q
\\
 \beta
&q\kappa_1-\kappa_2
\end{pmatrix}
\\
 \Delta&=&(q\kappa_1-\kappa_2)(q\theta_1-\theta_2)+q\beta
\end{eqnarray*}
A computation using (\ref{atilde}) gives $w\tilde y=q$, and 
$$w=R/\Delta,\qquad \tilde y=\frac{q\Delta}{R}$$ with 
\begin{eqnarray*}
\xi y^2\Delta
&=&\xi y^2(q\theta_1-\theta_2)(q\kappa_1-\kappa_2)+q\xi y^2\beta
\\
&=&\xi y^2(q\theta_1-\theta_2)(q\kappa_1-\kappa_2)
\\&&+qy\xi (-\kappa_1\kappa_2\sigma_3-\theta_1b-\theta_2a)\quad \text{using(\ref{beta})}
\\
&=& -q\kappa_1\kappa_2\theta_1\xi^2(y-c_3)(y-c_4)
\\&&+\xi ((q^2\kappa_1\theta_1+\kappa_2\theta_2)y^2-q\kappa_1\kappa_2\sigma_3y+2q\theta_1\theta_2)\\
&&-q\theta_2(y-c_1)(y-c_2)\quad \text{using (\ref{a}), (\ref{b})},
\end{eqnarray*}

and
\begin{eqnarray*}
\xi yR&=&\xi y((q\kappa_1-\kappa_2)(-aq^2+bq-q(q\kappa_1-\kappa_2)y)-q^2\gamma)\\
&=&\xi y(-q^3\kappa_1a-q\kappa_2b-q(q\kappa_1-\kappa_2)^2y-q^2\kappa_1\kappa_2\sigma_1)\quad\text{by (\ref{gamma})}\\
&=&q^3\kappa_1(\kappa_1\xi y^2-(y-c_1)(y-c_2)+\theta_1\xi)\\&&
+q\kappa_2(\kappa_2\xi y^2-\kappa_1\kappa_2\xi^2(y-c_3)(y-c_4)+\theta_2\xi)\\
&& -q(q\kappa_1-\kappa_2)^2y^2\xi-q^2\kappa_1\kappa_2\sigma_1\xi y\quad\text{by (\ref{a}), (\ref{b})}\\
&=&-q\kappa_1\kappa_2^2\xi^2(y-c_3)(y-c_4)
\\&&+\xi(2q^2\kappa_1\kappa_2y^2-q^2\kappa_1\kappa_2\sigma_1y
+q(q^2\kappa_1\theta_1+\kappa_2\theta_2))\\
&&-q^3\kappa_1(y-c_1)(y-c_2).
\end{eqnarray*}
Therefore 
\begin{eqnarray}\label{ytilde}
\tilde y&=&\frac{S}{yT}
\\ \nonumber\text{with}\\
S&=& -q\kappa_1\kappa_2\theta_1\xi^2(y-c_3)(y-c_4)\nonumber
\\&&+\xi ((q^2\kappa_1\theta_1+\kappa_2\theta_2)y^2-q\kappa_1\kappa_2\sigma_3y+2q\theta_1\theta_2)\nonumber\\
&&-q\theta_2(y-c_1)(y-c_2)\nonumber
\\
T&=&-\kappa_1\kappa_2^2\xi^2(y-c_3)(y-c_4)\nonumber
\\&&+\xi(2q\kappa_1\kappa_2y^2-q\kappa_1\kappa_2\sigma_1y
+(q^2\kappa_1\theta_1+\kappa_2\theta_2))\nonumber\\
&&-q^2\kappa_1(y-c_1)(y-c_2).\nonumber
\end{eqnarray}

In order to find $\tilde \xi$ we proceed as follows. From the definition one has
$$\tilde \xi=\frac{(\tilde y-c_1)(\tilde y-c_2)}{q\kappa_1 \tilde y^2+\tilde a\tilde y+q\theta_1}$$
where $\tilde a$ can be computed from (\ref{atilde}). Plugging in the values of 
$ a,b,\beta,\gamma$, and $\tilde y$ given by (\ref{a}), (\ref{b}), (\ref{beta}) (\ref{gamma}), (\ref{ytilde}) gives an unwieldy expression as a rational fraction in  the variables $y,\xi$. One can simplify this expression using the following observation. Computing the determinant 
$\text{det}\, \tilde A(\tilde y)=q\,\text{det}\, A(\tilde y)$ from (\ref{atilde}) gives the identity
\begin{equation}\label{det}
q\kappa_1\kappa_2(\tilde y-c_1)(\tilde y-c_2)(\tilde y-c_3)(\tilde y-c_4)=
(q\kappa_1 \tilde y^2+\tilde a\tilde y+q\theta_1)(\kappa_2 \tilde y^2+\tilde b\tilde y+\theta_2)
\end{equation}
where  $\tilde y,\tilde a,\tilde b$ are rational expressions in $y,\xi$. Of course this identity is compatible with the alternative expression 
$$\tilde \xi=\frac{\kappa_2 \tilde y^2+\tilde b\tilde y+\theta_2}{q\kappa_1\kappa_2(\tilde y-c_3)(\tilde y-c_4)}.$$
 Replacing $\tilde y$, $\tilde a$, $\tilde b$  in (\ref{det}) by their values and reducing the denominators, we get a polynomial identity between the numerators. This leads us to expect that the numerator of each expression $\tilde y-c_i$ factorizes into a product ot two terms, which divide respectively the numerators of 
$q\kappa_1 \tilde y^2+\tilde a\tilde y+q\theta_1$ and $\kappa_2 \tilde y^2+\tilde b\tilde y+\theta_2$. Since the degree of the numerator of $\tilde y-c_i$ in the variable $\xi$ is two, we expect a factorization into two terms of degree 1 in $\xi$. Indeed a direct computation  confirms that
\begin{eqnarray*}
S-c_1yT&=&\kappa_1\kappa_2\xi^2(c_1y\kappa_2-q\theta_1)(y-c_3)(y-c_4)\\&&
+\xi \left[-2qc_1\kappa_1\kappa_2y^3+(q^2\kappa_1\theta_1+\kappa_2\theta_2
+qc_1\kappa_1\kappa_2\sigma_1)y^2\right]\\&&
+\xi\left[(-q\kappa_1\kappa_2\sigma_3+c_1(q^2\kappa_1\theta_1+\kappa_2\theta_2))y+2q\theta_1\theta_2\right]\\&&
+q(qc_1\kappa_1y-\theta_2)(y-c_1)(y-c_2)\\&=&
(\xi(c_1y\kappa_2-q\theta_1)-qc_1(y-c_2))\times\\&&\times
(\kappa_1\kappa_2\xi(y-c_3)(y-c_4)-\frac{1}{c_1}(qc_1\kappa_1y-\theta_2)(y-c_1))
\end{eqnarray*}
and similar expressions for the other $c_i$, e.g. 

\begin{eqnarray*}
S-c_3yT&=&\kappa_1\kappa_2\xi^2(c_3y\kappa_2-q\theta_1)(y-c_3)(y-c_4)\\&&
+\xi \left[-2qc_3\kappa_1\kappa_2y^3+(q^2\kappa_1\theta_1+\kappa_2\theta_2
+qc_3\kappa_1\kappa_2\sigma_1)y^2\right]\\&&
+\xi\left[(-q\kappa_1\kappa_2\sigma_3+c_3(q^2\kappa_1\theta_1+\kappa_2\theta_2))y+2q\theta_1\theta_2\right]\\&&
+q(qc_3\kappa_1y-\theta_2)(y-c_1)(y-c_2)\\&=&
(\xi(y-c_4)-\frac{1}{\kappa_1\kappa_2 c_3}(qc_3\kappa_1y-\theta_2))\times\\&&\times
(\kappa_1\kappa_2\xi(y-c_3)(c_3y\kappa_2-q\theta_1)-\kappa_1\kappa_2 c_3(y-c_1)(y-c_2))
\end{eqnarray*} 
One checks then that one of these  factors divides the numerator of $q\kappa_1 \tilde y^2+\tilde a\tilde y+q\theta_1$.
This yields the  simpler expression 
\begin{equation}\label{xitilde}\tilde \xi=\frac{c_1c_2}{\kappa_1\theta_1\xi}
\frac{(\xi(y-\frac{q\theta_1}{c_1\kappa_2})-\frac{q}{\kappa_2}(y-c_2))(\xi(y-\frac{q\theta_1}{c_2\kappa_2})-\frac{q}{\kappa_2}(y-c_1))}
{(\xi(y-c_4)-\frac{q}{\kappa_2}(y-\frac{\theta_2}{qc_3\kappa_1}))
(\xi(y-c_3)-\frac{q}{\kappa_2}(y-\frac{\theta_2}{qc_4\kappa_1}))}.\end{equation}

\begin{proposition}\label{birat} Formulas  (\ref{ytilde}) and (\ref{xitilde}) express $\phi$ as a birational transformation $\phi:(y,\xi)\to(\tilde y,\tilde \xi)$.
\end{proposition}
\subsection{Blowing up the spaces}
The map $\phi$, whose expression as a birational transformation we obtained in the preceding section, can be extended to a regular map between two surfaces,  by blowing up the spaces 
$X_{\kappa_1,\kappa_2,\theta_1,\theta_2}$ and $X_{q\kappa_1,\kappa_2,q\theta_1,\theta_2}$.
For this we first extend the map to the  
 set  of endomorphisms such that
 $E_2= F_1$. As noted before, this set  corresponds to the blown up point
with coordinates
$(0,c_1c_2/\theta_2)$.
Take a matrix in a neighbourhood of  $E_2=F_1$ in the form
$$\begin{pmatrix}\kappa_1z^2+\bar az+\theta_2 & z-y\\
z(\bar \gamma z+\bar \beta)&\kappa_2z^2+\bar bz+\theta_1 
\end{pmatrix}$$
with $y$ small  and $\xi-c_1c_2\kappa_1/\theta_2=\lambda y=(y-c_1)(y-c_2)/z_1$, or
$$\bar a=-\kappa_1y+\frac{\kappa_1 y-\kappa_1(c_1+c_2)-\lambda\theta_2}{\lambda y+c_1c_2\kappa_1/\theta_2}$$ 
We can conjugate the matrix  to put it in the form 
$$\begin{pmatrix}\kappa_1z^2+az+\theta_1 & z-y\\
z(\gamma z+\beta)&\kappa_2z^2+bz+\theta_2 
\end{pmatrix}$$
with 
\begin{eqnarray*}
a&=&\bar a+\frac{\theta_2-\theta_1}{y}
\\ b&=&\bar b+\frac{\theta_1-\theta_2}{y}
\\ \gamma &=&\sigma_1+\kappa_1\bar b+\kappa_2\bar a+\frac{(\kappa_1-\kappa_2)(\theta_1-\theta_2)}{y}
\\ \beta &=&-\frac{\sigma_3+\theta_1\bar b+\theta_2\bar a}{y}
-\frac{(\theta_1-\theta_2)^2}{y^2}
\end{eqnarray*}
then putting  
$$\Delta=(q\theta_1-\theta_2)(q\kappa_1-\kappa_2)+q\beta$$
$$\tilde y=\frac{\Delta}{(q\kappa_1-\kappa_2)(-qa-y(q\kappa_1-\kappa_2))+(-q\gamma+
(q\kappa_1-\kappa_2)b)}$$
we obtain
$$\tilde a=\frac{(q\kappa_1-\kappa_2)(q\theta_1-\theta_2)qa-qy\beta(q\kappa_1-\kappa_2)+q\gamma(q\theta_1-\theta_2)+q\beta b}{\Delta}$$
We keep only the main terms as $y\to 0$ to get
$$\begin{array}{rcl}\Delta&\sim&-q\frac{(\theta_1-\theta_2)^2}{y^2}\\
\tilde y&\sim&-q\frac{(\theta_1-\theta_2)}{(q^2\kappa_1-\kappa_2)y}
\\ \tilde a&\sim&\frac{\theta_1-\theta_2}{y}\end{array}$$
therefore
$$\kappa_1\tilde z_1\sim q\kappa_1\tilde y^2+\tilde y\tilde a\sim \kappa_2\tilde y\frac{(\theta_1-\theta_2)}{(q^2\kappa_1-\kappa_2)y}\sim \kappa_2\tilde y^2/q$$
or $$\tilde \xi=q/\kappa_2$$
Thus we see that the map can be extended if we blow up the final space 
at $(\infty,q/\kappa_2)$. A similar computation shows that the inverse map also extends to the curve $\tilde E_2=\tilde F_1$ if we blow up the initial space at the point $(\infty,q/\kappa_2)$. 
\begin{proposition}
Let $Y_{\kappa_1,\kappa_2,\theta_1,\theta_2}$ denote 
${\bf P}^1\times {\bf P}^1$ blown up at the eight points  $$
\begin{matrix}(c_1,0),\quad (c_2,0),\quad(c_3,\infty),\quad(c_4,\infty)\\
(0,c_1c_2/\theta_1),\quad(0,c_1c_2/\theta_2),\quad(\infty,1/\kappa_1), 
\quad(\infty,q/\kappa_2)\end{matrix}$$
minus the strict transforms of the four lines $y=0,\infty$, $\xi=0,\infty$.
then the map $\phi$ extends to a regular isomorphism between the spaces 
$Y_{\kappa_1,\kappa_2,\theta_1,\theta_2}$ and 
$Y_{q\kappa_1,\kappa_2,q\theta_1,\theta_2}$.

\end{proposition}

$$\resizebox{7 cm}{!}{\begin{picture}(0,0)%
\includegraphics{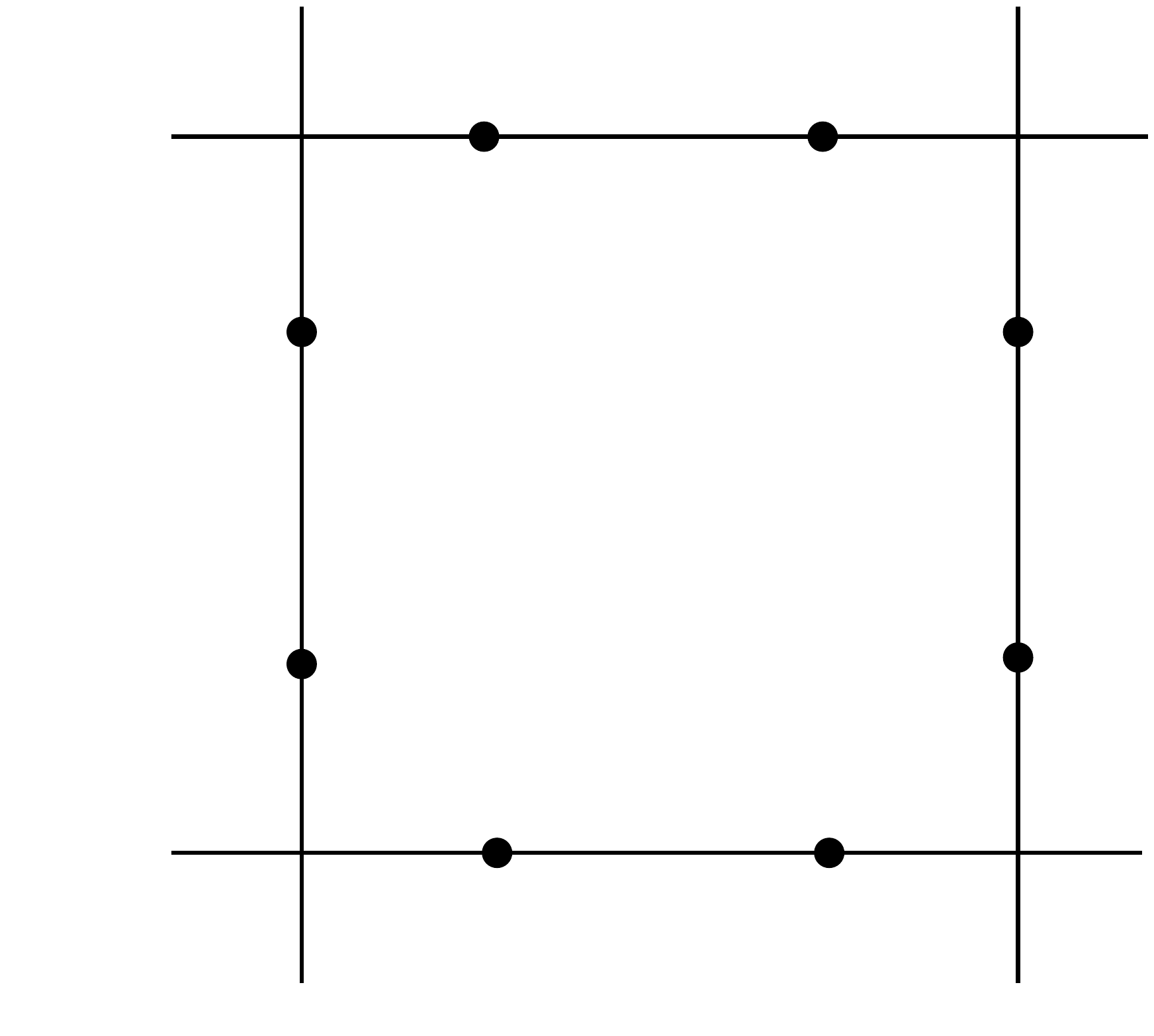}%
\end{picture}%
\setlength{\unitlength}{4144sp}%
\begingroup\makeatletter\ifx\SetFigFontNFSS\undefined%
\gdef\SetFigFontNFSS#1#2#3#4#5{%
  \reset@font\fontsize{#1}{#2pt}%
  \fontfamily{#3}\fontseries{#4}\fontshape{#5}%
  \selectfont}%
\fi\endgroup%
\begin{picture}(7968,7149)(1066,-7177)
\put(4006,-6496){\makebox(0,0)[lb]{\smash{{\SetFigFontNFSS{20}{24.0}{\rmdefault}{\mddefault}{\updefault}{\color[rgb]{0,0,0}$(c_1,0)$}%
}}}}
\put(6256,-6541){\makebox(0,0)[lb]{\smash{{\SetFigFontNFSS{20}{24.0}{\rmdefault}{\mddefault}{\updefault}{\color[rgb]{0,0,0}$(c_2,0)$}%
}}}}
\put(3646,-646){\makebox(0,0)[lb]{\smash{{\SetFigFontNFSS{20}{24.0}{\rmdefault}{\mddefault}{\updefault}{\color[rgb]{0,0,0}$(c_3,\infty)$}%
}}}}
\put(6031,-646){\makebox(0,0)[lb]{\smash{{\SetFigFontNFSS{20}{24.0}{\rmdefault}{\mddefault}{\updefault}{\color[rgb]{0,0,0}$(c_4,\infty)$}%
}}}}
\put(1081,-2446){\makebox(0,0)[lb]{\smash{{\SetFigFontNFSS{20}{24.0}{\rmdefault}{\mddefault}{\updefault}{\color[rgb]{0,0,0}$(0,c_1c_2/\theta_1)$}%
}}}}
\put(1126,-4696){\makebox(0,0)[lb]{\smash{{\SetFigFontNFSS{20}{24.0}{\rmdefault}{\mddefault}{\updefault}{\color[rgb]{0,0,0}$(0,c_1c_2/\theta_2)$}%
}}}}
\put(8326,-2401){\makebox(0,0)[lb]{\smash{{\SetFigFontNFSS{20}{24.0}{\rmdefault}{\mddefault}{\updefault}{\color[rgb]{0,0,0}$(\infty,1/\kappa_1)$}%
}}}}
\put(8281,-4651){\makebox(0,0)[lb]{\smash{{\SetFigFontNFSS{20}{24.0}{\rmdefault}{\mddefault}{\updefault}{\color[rgb]{0,0,0}$(\infty,q/\kappa_2)$}%
}}}}
\put(2701,-7036){\makebox(0,0)[lb]{\smash{{\SetFigFontNFSS{20}{24.0}{\rmdefault}{\mddefault}{\updefault}{\color[rgb]{0,0,0}$y=0$}%
}}}}
\put(1351,-6136){\makebox(0,0)[lb]{\smash{{\SetFigFontNFSS{20}{24.0}{\rmdefault}{\mddefault}{\updefault}{\color[rgb]{0,0,0}$\xi=0$}%
}}}}
\put(7426,-6991){\makebox(0,0)[lb]{\smash{{\SetFigFontNFSS{20}{24.0}{\rmdefault}{\mddefault}{\updefault}{\color[rgb]{0,0,0}$y=\infty$}%
}}}}
\put(1126,-1096){\makebox(0,0)[lb]{\smash{{\SetFigFontNFSS{20}{24.0}{\rmdefault}{\mddefault}{\updefault}{\color[rgb]{0,0,0}$\xi=\infty$}%
}}}}
\end{picture}%
} $$

\centerline{Fig. 3. The space $Y_{\kappa_1\kappa_2,\theta_1\theta_2}$}

\section{Discrete Painlev\'e equations}
\subsection{Sakai's theory}
We recall here the basic facts about the relation between surface theory and discrete Painlev\'e equations, according to Sakai \cite{S}.
One starts with a complex surface $X$, obtained by blowing up ${\bf P}^2({\bf C})$ at nine points 
$x_1,\ldots, x_9$,  thus $X$ comes with  a morphism $\pi:X\to{\bf P}^2({\bf C})$. Degenerate configurations, where some of the points $x_i$ are infinitely close, are allowed. The Picard group of such a surface is a free abelian group of  rank $10$,  generated by $\E_0=\pi^{-1}(L)$, the inverse image of a generic line
 in ${\bf P}^2({\bf C})$,  and 
 $\E_i=\pi^{-1}(x_i),i=1,\ldots,9$. The intersection form on Pic($X$) satisfies
$$ \E_0.\E_0=1,\quad  \E_i.\E_i=-1,\ i=1,\ldots,9,\quad  \E_i.\E_j=0,\ i\ne j$$

The canonical divisor class is
$$\K_X=-\delta=-3\E_0+\E_1+\ldots+\E_9$$
and its orthogonal $\delta^{\perp}$ in Pic($X$) is a rank 9 lattice, generated by
$$\E_i-\E_{i+1},i=1,\ldots,8, \quad \E_0-\E_1-\E_2-\E_3$$
which form the simple roots in an affine root system of type $E_8^{(1)}$
according to the Dynkin diagram
 $${\tt    \setlength{\unitlength}{0.60pt}
\begin{picture}(320,180)
\thinlines    \put(-80,20){ Fig. 1  Dynkin diagram of the affine root system $E_8^{(1)}$ inside Pic($X$)}

\put(-144,55){{$\E_1-\E_2$}}
\put(-64,55){{$\E_2-\E_3$}}
\put(16,55){{$\E_3-\E_4$}}
\put(96,55){{$\E_4-\E_5$}}
\put(176,55){{$\E_5-\E_6$}}
\put(256,55){{$\E_6-\E_7$}}
\put(336,55){{$\E_7-\E_8$}}
\put(416,55){{$\E_8-\E_9$}}

\put(-115,83){\line(1,0){70}}
 \put(-35,83){\line(1,0){70}}
\put(45,83){\line(1,0){70}}
\put(125,83){\line(1,0){70}}
\put(205,83){\line(1,0){70}}
\put(285,83){\line(1,0){70}}
\put(365,83){\line(1,0){70}}

 \put(-119,83){\circle{8}}
\put(-39,83){\circle{8}}
\put(41,83){\circle{8}}
 \put(121,83){\circle{8}}
\put(201,83){\circle{8}}
  \put(281,83){\circle{8}}
    \put(361,83){\circle{8}}
\put(441,83){\circle{8}}
              
\put(41,87){\line(0,1){60}}
\put(41,151){\circle{8}}
\put(-20,165){{$\E_0-\E_1-\E_2-\E_3$}}

              \end{picture}}
$$

The surface $X$ is called a generalized Halphen surface by Sakai, if the linear system $|-\K_X|$ has a unique divisor, whose irreducible components are orthogonal to $\K_X$. Let $\D_i,i\in I$ be these components, then they generate an affine subroot system $R$ in 
$\delta^{\perp}\subset\text{Pic}(X)$, such that its orthogonal is also an affine subroot system $R^{\perp}\subset\delta^{\perp}$. These surfaces are classified according to their corresponding root systems.
Let $X$ and $\tilde X$ be two such surfaces, corresponding to the same root system
$R$, where the irreducible components of the canonical divisor are labelled
$\D_i, i\in I$, resp. $\tilde{\D}_i,i\in I$, 
then 
Sakai's main result is that (up to some little twists, see \cite{S}, Theorem 25 for a precise statement),
 every isometry of the Picard groups $\phi:\text{Pic}(X)\to\text{Pic}(\tilde X)$ which preserves the components of the anticanonical divisor,
$\phi((\D_i,i\in I))=(\tilde{\D}_i,i\in I)$ comes from a unique isomorphism of surfaces $\phi:X\to\tilde X$. Such isometries are given by elements of  
the extended Weyl group of the dual root system $R^{\perp}$ (i.e. the group generated by the Weyl group of $R^{\perp}$ and the automorphisms of the Dynkin diagram). When viewed at the level of ${\bf P}^2({\bf C})$, the
isomorphisms $\phi:X\to\tilde X$ give rise to birational transformations ${\bf P}^2({\bf C})\to{\bf P}^2({\bf C})$. The discrete Painlev\'e equations correspond to  the translations in the affine Weyl group $W(R^{\perp})$.
On the other hand, transformations in the extended Weyl group of $R$ give 
B\"acklund transformations.
As an example, Sakai treats the case of the birational tranformation 
$$\rho:(f,g)\to(\bar f,\bar g)$$  determined by
\begin{eqnarray*}\bar g&=&\frac{b_3b_4}{g}\frac{(f-b_5)(f-b_6)}{(f-b_7)(f-b_8)}\\
\bar f&=&\frac{b_7b_8}{f}\frac{(\bar g-qb_1)(\bar g-qb_2)}{(\bar g-b_3)(\bar g-b_4)}
\end{eqnarray*}
where $b_1,\ldots,b_8$ are some complex parameters and $$q=\frac{b_3b_4b_5b_6}{b_1b_2b_7b_8}.$$
He shows that this transformation corresponds to an isomorphism between $A_3^{(1)}$ surfaces, which arises from a translation in an affine Weyl group $D_5^{(1)}$. Here $A_3^{(1)}$ and $D_5^{(1)}$
are two orthogonal affine root systems in $E_8^{(1)}$.
Remark that this transformation is actually considered as a birational tranformation of 
${\bf P}^1({\bf C})\times{\bf P}^1({\bf C})$. Since ${\bf P}^2({\bf C})$ blown up at 9 points is isomorphic to 
${\bf P}^1({\bf C})\times{\bf P}^1({\bf C})$ blown up at 8 points the two pictures are equivalent, but it is sometimes easier to do formal reasoning on ${\bf P}^2({\bf C})$ and to perform actual computations in 
${\bf P}^1({\bf C})\times{\bf P}^1({\bf C})$.

As explained in \cite{JS}, the transformation $\rho$ can be considered as a step in a nonlinear difference equation
$$(f_{n+1},g_{n+1})=\rho(f_n,g_n)$$
and this difference equation leads, by taking a suitable continuous limit, to a rational differential system
\begin{eqnarray*}
df/dx&=&p(x,f,g)\\
dg/dx&=&q(x,f,g)
\end{eqnarray*}
which is equivalent to a Hamiltonian version of the Painlev\'e VI equation.
The degenerescence scheme of discrete Painlev\'e equations corresponds to the inclusion of the root systems in $E_8^{(1)}$. The degenerescence to Painlev\'e differential equations is also discussed by Sakai.

\subsection{}According to the classification of Sakai, the spaces
$Y_{\kappa_1,\kappa_2,\theta_1,\theta_2}$ are examples of $A_3^{(1)}$ surfaces. We will now identify the map $\phi$ between  $Y_{\kappa_1,\kappa_2,\theta_1,\theta_2}$ and 
$Y_{q\kappa_1,\kappa_2,q\theta_1,\theta_2}$ we have constructed as one of the discrete Painlev\'e maps. For this we will use intersection theory for surfaces obtained as blow-ups.
Denote by 
$\mathcal F_i,i=1,\ldots ,8$ the divisors which are inverse images, respectively,
 of the points $$\begin{matrix}(c_1,0),\quad(c_2,0),\quad(c_3,\infty),\quad(c_4,\infty),\\ (0,c_1c_2/\theta_1),\quad(0,c_1c_2/\theta_2),\quad(\infty,1/\kappa_1),
\quad(\infty,q/\kappa_2)\end{matrix}$$ in the blow ups.
Similarly we denote by $\tilde{\mathcal F}_i,i=1,\ldots ,8$ the ones in the 
space $Y_{q\kappa_1,\kappa_2,q\theta_1,\theta_2}$. 

It will be convenient to identify, as in Sakai \cite{S}, the surfaces with   ${\bf P}^2$ blown up at   9 points, called $x_i$.
This identification is as follows. Consider 
 ${\bf P}^1\times {\bf P}^1$ with coordinates denoted by  $(f,g)$, blown up at the eight points 
$$\begin{matrix}(0,b_1),\quad(0,b_2),\quad(\infty,b_3),\quad(\infty,b_4),\\ (b_5,0),\quad(b_6,0),\quad(b_7,\infty),
\quad(b_8,\infty)\end{matrix}$$
We  make  the explicit identification
\begin{equation}
\begin{matrix}\label{corr}
c_1=b_1\\c_2=b_2\\c_3=b_3\\c_4=b_4\end{matrix}\qquad\qquad
\begin{matrix}c_1c_2/\theta_1=b_5\\
c_1c_2/\theta_2=b_6\\1/\kappa_1=b_7\\
q/\kappa_2=b_8\end{matrix}\qquad f=\xi\qquad g=y
\end{equation}between our space and Sakai's. Consider also
 ${\bf P}^2$ blown up at the eight points
$$\begin{matrix}x_1=(a_1:0:1),\quad x_2=(a_2:0:1),\quad x_3=(a_3:0:1),\quad x_4=(0:1:a_4),\\ x_5=(0:1:a_5),\quad x_6=(1:a_6:0),\quad x_7=(1:a_7:0),\quad x_8=(0:1:0)
\end{matrix}$$
further blown up on $\pi^{-1}(0:1:0)$ at $$x_9=\left(\frac{x}{z},\frac{z}{x}\right)=(0,a_8)$$
The two spaces are identified according to the transformation
\begin{equation}\label{p1p2}
\left(\begin{matrix}b_1&b_2&b_3&b_4
 \\b_5&b_6&b_7&b_8\end{matrix}\ ;f,g\right)=
\left(\begin{matrix}1/a_2&1/a_3&1/a_1&a_8
 \\a_6&a_7&-1/a_4a_1&-1/a_5a_1\end{matrix}\ ;\frac{y}{x-a_1z},\frac{z}{x}\right)
\end{equation}
See pages 173 and 215 in \cite{S}.
The correspondance between divisors is

\begin{eqnarray*}\E_2=\F_1,\ \E_3=\F_2,\ \E_4=\F_5,\ \E_5=\F_6,\ \E_6=\F_7,\\
 \E_7=\F_8,\ \E_9=\F_4,\ \E_0-\E_1-\E_8=\F_3
\end{eqnarray*} $\E_1$ is the line  $y=a_3$, and $\E_8-\E_9$ the line $\xi=\infty$.

The components of the anticanonical cycle are 

\begin{equation}\label{D}\begin{array}{rclcrcl}\D_0&=&\E_8-\E_9&\qquad&
\D_1&=&\E_0-\E_6-\E_7-\E_8 \\ \D_2&=&\E_0-\E_1-\E_2-\E_3 &\qquad&
\D_3&=&\E_0-\E_4-\E_5-\E_8
\end{array}
\end{equation}

which form a root system of type $A_3^{(1)}$, with Dynkin diagram
$${\tt    \setlength{\unitlength}{0.60pt}
\begin{picture}(263,120)
\thinlines         
 
\put(-5,15){$\D_1$}
\put(0,40){\circle{8}}
\put(5,40){\line(1,0){80}}
\put(4,42){\line(2,1){81}}
\put(90,40){\circle{8}}
\put(85,15){$\D_2$}
 \put(95,40){\line(1,0){80}}
\put(90,82){\circle{8}}
\put(85,100){$\D_0$} 
\put(176,42){\line(-2,1){81}}       
   \put(180,40){\circle{8}}  
\put(175,15){$\D_3$}         
\end{picture}}
$$

The orthogonal is a root system of type $D_5^{(1)}$, with generators
\begin{equation}\label{alpha}
\begin{array}{rclcrclcrcl}\alpha_0&=&\E_0-\E_1-\E_8-\E_9,&\quad& \alpha_1&=&\E_2-\E_3,&\quad& \alpha_2&=&\E_1-\E_2\\
\alpha_3&=&\E_0-\E_1-\E_4-\E_6,&\quad&  \alpha_4&=&\E_6-\E_7,&\quad &\alpha_5&=&\E_4-\E_5\end{array}
\end{equation}
forming the Dynkin diagram
$${\tt    \setlength{\unitlength}{0.60pt}
\begin{picture}(263,150)
\thinlines    
\put(57,83){\line(-1,2){20}}
\put(0,124){{$\alpha_0$}}
\put(34,126){\circle{8}}
\put(57,77){\line(-1,-2){20}}
\put(0,34){{$\alpha_1$}}
\put(34,34){\circle{8}}
\put(60,80){\circle{8}}
\put(25,80){{$\alpha_2$}}
\put(65,80){\line(1,0){80}}
\put(185,124){{$\alpha_4$}}
\put(150,80){\circle{8}}
\put(165,80){{$\alpha_3$}}
\put(153,83){\line(1,2){20}}
    \put(153,77){\line(1,-2){20}} 
\put(176,126){\circle{8}} 
\put(176,34){\circle{8}}      
\put(185,34){{$\alpha_5$}}       
\end{picture}}
$$
We denote by $\tilde \E$ and $\tilde\F$ the objects corresponding to the target space
of the map $\phi$.
We saw that $\phi$ maps $\F_6$  to $\tilde\F_8$ and $\F_7$ to $\tilde \F_5$.
We will compute the matrix of the map induced by  $\phi$ between the Picard groups.
\subsection{}
Let us start by computing the image of the curve  $\F_1=\E_2$. This curve consists of matrices of the form 
$$A(z)=\begin{pmatrix}
\kappa_1z^2+az+\theta_1& w(z-c_1)\\w^{-1}z(\gamma z+\beta)&
\kappa_2(z-c_1)(z-\theta_2/c_1\kappa_2)
\end{pmatrix}$$
where  $a$ is a  local  parameter along the curve and $w$ an irrelevant constant which can be removed by conjugation with a diagonal matrix but is convenient for computations (this is the "gauge freedom" in \cite{JS}). The parameters
 $\gamma$ and $\beta$ are 
 affine functions of $a$. Recall that
$$\Delta=(q \kappa_1-\kappa_2)(q\theta_1-\theta_2)+q\beta$$
then
$$\tilde y=\frac{\Delta}{qa(q\kappa_1-\kappa_2)-(q\kappa_1-\kappa_2)(-c_1-\theta_1/a_1\kappa_1)+c_1(q\kappa_1-\kappa_2)^2
+\gamma(q\kappa_1-\kappa_2)}$$
The  critical value to be investigated are 
$\tilde y=0,c_i,\infty$.
The  case $\tilde y=0$ corresponds to  $\tilde E_1=\tilde F_2$ so there is an intersection with  
$\tilde\F_6=\E_5$. In the case  $\tilde y=c_1$, the matrix $A(c_1)$ has the form
$\begin{pmatrix}
S&0 \\T&
0
\end{pmatrix}$, therefore $$\tilde A(c_1)=\begin{pmatrix}
(q\theta_1-\theta_2)W&-qW \\(q\theta_1-\theta_2)Z&-qZ
\end{pmatrix}$$
with  $W=(q\kappa_1-\kappa_2)qS+\beta T$, thus $\tilde y=c_1$ if $W=0$, and then
  $\tilde z_1=0$, therefore the image does not intersect
 $\tilde \F_1=\tilde \E_2$.
In the case  $y=c_i$, with $i=1,2,3$, the matrix $A(c_i)$ has rank 1 so
$$A(c_i)=\begin{pmatrix}
rG& G\\rH&
H
\end{pmatrix}$$
for some constants  $H$ and $G$  (independent of $a$), and
$$\tilde A(c_i)=\begin{pmatrix}
((q\theta_1-\theta_2)r-\beta)((q\kappa_1-\kappa_2)G-qH)& (qr+q\kappa_1-\kappa_2)((q\kappa_1-\kappa_2)G-qH)\\
(q\theta_1-\theta_2)r-\beta)(\beta G+q\theta_1-\theta_2)&(qr+q\kappa_1-\kappa_2)(\beta G+q\theta_1-\theta_2)
\end{pmatrix}$$
The $12$ coefficient  vanishes only if $qr+q\kappa_1-\kappa_2=0$, and then $\tilde z_2=0$, therefore the image intersects 
$\tilde \F_i$ if $i=2$ but not $i=3,4$, and intersects
  $\tilde\E_1$. Finally for some value of  $a$ the equation
 $\tilde y=\infty$ holds and the image intersects
$\tilde \F_7$.
We conclude from these considerations that, in the Picard groups of the surfaces,
$$\phi(\E_2)=2\tilde\E_0-\tilde\E_1-\tilde\E_3-\tilde\E_5-\tilde\E_6-\tilde\E_8.$$

Similar tedious computations,  which we omit, lead  
 finally to the matrix induced by  the map $\phi$ between the Picard groups, which is, from the basis
$(\E_i)_{i=0,\ldots,9}$ to $(\tilde \E_i)_{i=0,\ldots,9}$ given by
$$\begin{pmatrix}6&2&2&2&3&0&0&3&2&1\\
-2&0&-1&-1&-1&0&0&-1&-1&0\\-2&-1&0&-1&-1&0&0&-1&-1&0\\
-2&-1&-1&0&-1&0&0&-1&-1&0\\0&0&0&0&0&0&1&0&0&0\\
-3&-1&-1&-1&-2&0&0&-1&-1&-1\\
-3&-1&-1&-1&-1&0&0&-2&-1&-1\\
0&0&0&0&0&1&0&0&0&0\\
-2&-1&-1&-1&-1&0&0&-1&0&0\\
-1&0&0&0&-1&0&0&-1&0&0
\end{pmatrix}$$

One checks that this transformation preserves the components $\D_i$ (cf (\ref{D}))
$$\begin{array}{rclcrcl}\phi(\D_0)&=&\tilde\D_2&\qquad&\phi(\D_2)&=&\tilde\D_0\\
\phi(\D_1)&=&\tilde\D_3&\qquad&\phi(\D_3)&=&\tilde\D_1
\end{array}
$$
while its action on the orthogonal root system (\ref{alpha}) is a translation in the affine Weyl group, namely

\begin{eqnarray*}\phi(\alpha_0)=\alpha_0,\qquad \phi(\alpha_1)=\alpha_1,\qquad \phi(\alpha_2)=\alpha_2,\qquad \\
\phi(\alpha_3)=\alpha_3,\qquad \phi(\alpha_4)=\alpha_4-\delta,\qquad
\phi(\alpha_5)=\alpha_5+\delta
\end{eqnarray*}
with
$$\delta=\alpha_0+\alpha_1+2\alpha_2+2\alpha_3+\alpha_4+\alpha_5$$

This  translation has the following reduced  decomposition, where the
 $w_i$ are the simple reflexions associated with the roots  $\alpha_i$,
and 
 $\sigma$ is the automorphism of the Dynkin  diagram such that
  $\sigma(\alpha_0)=\alpha_1,
 \sigma(\alpha_4)=\alpha_5$.
\begin{equation}\label{red}\phi=\sigma w_4w_3w_2w_0w_1w_2w_3w_4\end{equation}
We summarize the results of this section in the following.
\begin{proposition}
The transformation $\phi$ induces, on the Picard group of the surfaces, a translation whose reduced decomposition is given by 
$$\phi=\sigma w_4w_3w_2w_0w_1w_2w_3w_4$$
\end{proposition}
\begin{rem}
The discrete Painlev\'e equation of \cite{JS} corresponds to the translation
$\psi$ given by
\begin{eqnarray*}\psi(\alpha_i)=\alpha_i, \qquad i=0,1,4,5 \\
 \psi(\alpha_2)=\alpha_2-\delta,\qquad
\psi(\alpha_3)=\alpha_3+\delta
\end{eqnarray*}
The translations $\phi$ and $\psi$ are not conjugate in the Weyl group.
\end{rem}

\subsection{}As a check, we will obtain formulas (\ref{ytilde}), (\ref{xitilde})
using (\ref{red}) and the explicit expressions for the birational transformations corresponding to the reflexions. 
After Sakai, we have the following formulas for the  transformations corresponding to the roots
$$\begin{array}{c}w_0:\left(\begin{matrix}a_8& a_1& a_2& a_3  \\a_4&a_5&a_6&a_7\end{matrix}\ ;x:y:z\right)\to \\
\left(\begin{matrix}1/a_1&1/a_8&a_2&a_3\\a_1a_4&a_5a_1&a_6a_8&a_7a_8\end{matrix}\ 
;x(x-a_1z):y(x-z/a_8):z(x-a_1z)\right)\end{array}$$

 $w_1,w_4,w_5$ are the  transpositions of parameters
$$w_1:(a_2,a_3);\quad w_4:(a_4,a_5),\quad w_5:(a_6,a_7)$$
Passing to   ${\bf P}^1\times {\bf P}^1$ with (\ref{p1p2}), 
 $w_2$ and $w_3$ are
$$\begin{array}{c}w_2:\left(\begin{matrix}b_1&b_2&b_3&b_4
 \\b_5&b_6&b_7&b_8\end{matrix}\ ;f,g\right)\to \\
\left(\begin{matrix}b_3&b_2&b_1&b_4
 \\b_5\frac{b_3}{b_1}&b_6\frac{b_3}{b_1}&b_7&b_8\end{matrix}\ ;f\frac{g-b_3}{g-b_1},g\right)\end{array}$$
$$\begin{array}{c}w_3:\left(\begin{matrix}b_1&b_2&b_3&b_4
 \\b_5&b_6&b_7&b_8\end{matrix}\ ;f,g\right)\to \\
\left(\begin{matrix}b_1&b_2&\frac{b_5}{b_7}b_3&\frac{b_5}{b_7}b_4
 \\b_7&b_6&b_5&b_8\end{matrix}\ ;f,\frac{b_5}{b_7}g\frac{f-b_7}{f-b_5}\right)\end{array}$$
The automorphism  $\sigma$ is the inversion
$$\left(\begin{matrix}b_1&b_2&b_3&b_4
 \\b_5&b_6&b_8&b_7\end{matrix}\ ;f,g\right)\to
\left(\begin{matrix}b_1&b_2&b_3&b_4
 \\b_5&b_6&b_8&b_7\end{matrix}\ ;b_7b_8/f,b_3b_4/g\right)
$$

We can now compute 
 $\phi$ induced by the translation (\ref{red}). One gets

$$\bar f=\frac{b_1b_2b_8}{b_5f}\frac
{(f(g-\frac{b_1b_8}{b_5})-b_8(g-b_1))
(f(g-\frac{b_2b_8}{b_5})-b_8(g-b_2)}
{(f(g-b_3)-b_8(g-\frac{b_3b_5}{b_8})
(f(g-b_4)-b_8(g-\frac{b_4b_5}{b_8})}
$$
$$\bar g=\frac{b_5b_7}{g}\frac{f-b_5}{f-b_8}\frac{
f\frac{g\frac{f-b_8}{f-b_5}-b_3}{g\frac{f-b_8}{f-b_5}-b_1\frac{b_8}{b_5}}\frac{g\frac{f-b_8}{f-b_5}-b_4}{g\frac{f-b_8}{f-b_5}-b_2\frac{b_8}{b_5}}-\frac{b_3b_4b_5^2}{b_1b_2b_8}}
{f\frac{g\frac{f-b_8}{f-b_5}-b_3}{g\frac{f-b_8}{f-b_5}-b_1\frac{b_8}{b_5}}\frac{g\frac{f-b_8}{f-b_5}-b_4}{g\frac{f-b_8}{f-b_5}-b_2\frac{b_8}{b_5}}-b_5}$$
This gives (\ref{ytilde}), (\ref{xitilde}), using the correspondance 
(\ref{corr})
and after simplifying the expression for $\bar g$, noting that $f-b_5$ divides the denominator, while $f-b_8$ divide the numerator.

\subsection{}
Finally the results of the preceding section allow one to derive relations for the 
Verblunsky coefficients. Indeed one has 
$$\text{det}\ A_n(z)=q^nV(z)W(z)$$ 
thus one can take $c_1=\bar a/q,\quad c_2=1/b \quad c_3=\bar b/q,\quad c_4=1/b$. From (\ref{theta_n}) one sees that the coordinate $y$
corresponds to $\frac{\alpha_{n}(\bar a-\bar bq^{n})}{\alpha_{n+1}(a-bq^{n+1})}$.
According to equations (\ref{omega_n}), (\ref{omega_n*}), one finds that
$\kappa_1=bq^{n+1},\quad \kappa_2=aq,\quad \theta_1=\bar bq^n, \quad\theta_2=\bar a$.
 We conclude that these expressions satisfy a discrete Painlev\'e equation.

\section{Recovering the  differential equation}

We shall now explain how to obtain a differential system
 as a limit 
of the discrete dynamical system given by 
(\ref{ytilde}), (\ref{xitilde}).
For this we introduce a complex variable $t$, and consider the parameters
$(t\kappa_1,\kappa_2,t\theta_1,\theta_2,c_1,c_2,c_3,c_4)$. Consider the
map $\phi$ of Proposition
(\ref{birat}) as a birational map  $(y(tq),\xi(tq))=\phi(y(t),\xi(t))$. Using the parametrisation
\begin{eqnarray*}
&q=1-\epsilon, \\
&q\kappa_1=t(1+\epsilon K_1),\ \kappa_2=1+\epsilon K_2\\
&q\theta_1=t(1+\epsilon \Theta_1),\
\theta_2=1+\epsilon \Theta_2,\\
&c_i=1+\epsilon C_i, i=1,\ldots, 4\\ &y=1+\epsilon u,\xi=v
\end{eqnarray*}
 Letting $\epsilon$ go to zero gives a differential system
in the variable $t$.
We omit the details of the lengthy computation, and give the final result:

\begin{eqnarray*}
t(1-t)\frac{du}{dt}&=&tv(u-C_3)(u-C_4)-t(t-1)v^{-1}(u-C_1)(u-C_2)\\
\\
t(1-t)\frac{dv}{dt}&=&v\left[2u+2t(K_2-\Theta_1)+C_1+C_2-K_1-\Theta_1
 \right]+\\&&
2u(t-1)+C_2+C_1
-t(C_3+C_4)+\\&&v^{-1}\left[2ut-t(1+K_2)+2(1+\Theta_2)-C_3-C_4
\right]
\end{eqnarray*}

Finally, we come back to the original problem about the Fourier transform and the scattering problem for phase shifts which are quotients of Gamma functions, like
(\ref{Gamma}). 
First note that in the case of the weight $w$ given by (\ref{w}), one has, from
 (\ref{omega_n}),  (\ref{omega_n*}), and (\ref{VV}), (\ref{WW}),
 
\begin{eqnarray*}
\kappa_1=bq, \theta_1=\bar b, \kappa_2=aq, \theta_2=\bar a\\
c_1=\bar a/q,c_2=\bar b/q,c_3=1/b,c_4=1/a
\end{eqnarray*}

In order to recover  (\ref{Gamma}) as a limit of the phase shift associated with the weight $w$ we make the substitution
$q\to 1-\epsilon$, $a\to 1+\epsilon a$,  $b\to 1+\epsilon b$ where $a$ and $b$ are real numbers, before letting $\epsilon \to 0$. The final system is then

\begin{eqnarray*}
t(1-t)\frac{du}{dt}&=&tv(u+b)(u+a)-t(t-1)v^{-1}(u-a+1)(u-b+1)\\
\\
t(1-t)\frac{dv}{dt}&=&v\left[2u+2t(a-b-1)+a-b+3
 \right]+\\&&
2u(t-1)+a+b+2
+t(a+b)+\\&&v^{-1}\left[(2u-a)t+2+3a+b
\right]
\end{eqnarray*}


\end{document}